\documentclass[11pt, twoside]{article}
\usepackage{latexsym}
\usepackage{amsmath}
\usepackage{amssymb}
\usepackage[all]{xy}
\usepackage{amsfonts}
\usepackage{verbatim}
\usepackage{amsthm}
\usepackage{bm}
\usepackage{mathrsfs}
\usepackage{epsfig}
\usepackage{xy}
\usepackage{array}
\usepackage{stmaryrd}
\usepackage{graphicx,color}
\usepackage{xcolor}
\usepackage{tikz}
\usepackage{xcolor}
\usepackage[colorlinks=true,linkcolor=blue,citecolor=blue]{hyperref}
\usetikzlibrary{arrows,calc}
\usepackage{etex}
\usepackage{mathdots}
\usepackage{float}
\usepackage{graphics}
\usepackage{pdflscape}
\usepackage{CJK}
\usepackage{anysize,hyperref}
\input xypic
\xyoption{all}
\usepackage{extarrows}
\usepackage[perpage,symbol]{footmisc}
\usepackage{setspace}
\usepackage{enumitem}
\setenumerate[1]{itemsep=0pt,partopsep=0pt,parsep=\parskip,topsep=5pt}
\setitemize[1]{itemsep=0pt,partopsep=0pt,parsep=\parskip,topsep=5pt}
\setdescription{itemsep=0pt,partopsep=0pt,parsep=\parskip,topsep=5pt}

\def\A{\mathcal{A}}

\def\T{\mathcal{T}}

\def\C{\mathscr{C}}

\def\dr{\ar@{->}[r]}

\def\X{\mathscr{X}}

\def\add{\mbox{add}}

\renewcommand{\diagram}[3]{\matrix (#1) [matrix of math nodes,row
  sep={#2},column sep={#3},text height=1.5ex,text depth=0.25ex]}
\usepackage{pgf,tikz}

  \usetikzlibrary{arrows}
  \usetikzlibrary{decorations.markings}
  \usetikzlibrary{shapes}
  \usetikzlibrary{matrix}
  \usepackage{wasysym}
\begin{document}
\baselineskip=15pt
\title{\Large{\bf The axioms for right $\bm{(n+2)}$-angulated categories\footnotetext{Jing He was supported by the Hunan Provincial Natural Science Foundation of China (Grant No. 2023JJ40217). }}}
\medskip
\author{Jing He and Jiangsha Li}

\date{}

\maketitle
\def\blue{\color{blue}}
\def\red{\color{red}}

\newtheorem{theorem}{Theorem}[section]
\newtheorem{lemma}[theorem]{Lemma}
\newtheorem{corollary}[theorem]{Corollary}
\newtheorem{proposition}[theorem]{Proposition}
\newtheorem{conjecture}{Conjecture}
\theoremstyle{definition}
\newtheorem{definition}[theorem]{Definition}
\newtheorem{question}[theorem]{Question}
\newtheorem{remark}[theorem]{Remark}
\newtheorem{remark*}[]{Remark}
\newtheorem{example}[theorem]{Example}
\newtheorem{example*}[]{Example}
\newtheorem{condition}[theorem]{Condition}
\newtheorem{condition*}[]{Condition}
\newtheorem{construction}[theorem]{Construction}
\newtheorem{construction*}[]{Construction}

\newtheorem{assumption}[theorem]{Assumption}
\newtheorem{assumption*}[]{Assumption}

\baselineskip=17pt
\parindent=0.5cm

\begin{abstract}
\baselineskip=16pt
Drawing inspiration from the works of Beligiannis--Marmaridis and Lin, we refine the axioms for a right $(n+2)$-angulated category and give some examples of such categories. Interestingly, we show that the morphism axiom for a right $(n+2)$-angulated category is actually redundant. Moreover, we prove that the higher ``octahedral axiom" is equivalent to the mapping cone axiom for a right
$(n+2)$-angulated category.
\\[0.5cm]
\medskip
\end{abstract}

\pagestyle{myheadings}
\markboth{\rightline {\scriptsize J. He and J. Li}}
         {\leftline{\scriptsize  The axioms for right $(n+2)$-angulated categories}}

\section{Introduction}
The notion of a triangulated category was introduced in the mid 1960's by Verdier
in his thesis \cite{V}. Having their origins in algebraic geometry and algebraic topology, triangulated categories have now become indispensable in many different areas of mathematics.
Assem, Beligiannis and Marmaridis \cite{ABM,BM} defined
the notion of right triangulated category. Informally, a right triangulated category is a triangulated category whose suspension functor is not necessarily an automorphism.
Let ${\rm mod}\Lambda$ be the category of finitely generated right
$\Lambda$-modules over an artin algebra $\Lambda$.
Beligiannis and Marmaridis \cite{BM} proved that
any covariantly finite subcategory $\X$ of ${\rm mod}\Lambda$ induces a right triangulated category
on the quotient category ${\rm mod}\Lambda/\X$.
Later, this result was extended by Beligiannis and Reiten \cite{BR} to a more general situation, namely, if $\A$ is an abelian category and $\X$ is contravariantly finite in $\A$,
then the quotient category $\A/\X$ is a right triangulated category.

Recently, Geiss, Keller and Oppermann \cite{GKO} introduced the notion of an $(n+2)$-angulated category, which is a higher dimensional analogue of a triangulated category.
We note that the case $n=1$ corresponds to a triangulated category.
A primary source of examples for $(n+2)$-angulated categories is $n$-cluster tilting subcategories of triangulated categories that are closed under the $n$th power of the shift functor.
Building on the concept of a right triangulated category, Lin \cite{L2} defined a right $(n+2)$-angulated category and explored those that arise from a covariantly finite subcategory.

In this paper, we refine the axioms for a right $(n+2)$-angulated category as defined by Lin.
Moreover, we also provide some examples of a right $(n+2)$-angulated category.
The morphism axiom of an $(n+2)$-angulated category states that a morphism between
the bases of two $(n+2)$-angles can be extended to a morphism of $(n+2)$-angles.
Arentz-Hansen, Bergh and Thaule \cite{ahbt} showed that the morphism axiom for an $(n+2)$-angulated category is redundant. Based on this idea, we prove that the conclusion holds true in a right $(n+2)$-angulated category.
Bergh and Thaule \cite{bt} showed that the higher ``octahedral axiom" is equivalent to the mapping cone axiom in an $(n+2)$-angulated category. We explain how this generalizes to right $(n+2)$-angulated categories. However, our proof method differs from that used in the case of the $(n+2)$-angulated category. Additionally, we present another equivalent characterization of the mapping cone axiom.

The paper is organized as follows:
In Section 2, we provide an overview of the definitions of the $(n+2)$-angulated category, the $n$-cokernel, and the special $n$-cokernel.
In Section 3, we refine the definition of a right $(n+2)$-angulated category.
In Section 4, we present some examples of right $(n+2)$-angulated categories.
In Section 5, we prove the redundancy of the morphism axiom in a right $(n+2)$-angulated category.
In Section 6, we give some new equivalent statements of the
higher mapping cone axiom.

\section{Preliminaries}
In this section, we recall
the axioms for $(n+2)$-angulated categories as described in \cite{GKO,bt}.
Let $\C$ be an additive category with an automorphism $\Sigma: \C\rightarrow\C$, and $n$ be an positive integer. A sequence of objects and morphisms in $\C$ of the form
$$A_0 \xrightarrow{~a_0~} A_1 \xrightarrow{~a_1~} A_2 \xrightarrow{~a_2~} \cdots \xrightarrow{~a_{n-1}~}A_n \xrightarrow{~a_{n}~} A_{n+1} \xrightarrow{~a_{n+1}~} \Sigma A_0$$
is called an $(n+2)$-$\Sigma$-$sequence$. Its {\em left rotation} is in the form of the following
$$A_1 \xrightarrow{~a_1~} A_2 \xrightarrow{~a_2~} A_3 \xrightarrow{~a_3~} \cdots \xrightarrow{~a_{n}~}A_{n+1} \xrightarrow{~a_{n+1}~} \Sigma A_{0} \xrightarrow{(-1)^n\Sigma a_{0}} \Sigma A_1.$$
A $morphism$ of $(n+2)$-$\Sigma$-sequences is a sequence $f=(f_0, f_1, \cdots, f_{n+1})$ in $\C$ such that the following diagram
$$\xymatrix{
A_0 \ar[r]^{a_0}\ar[d]^{f_0} & A_1 \ar[r]^{a_1}\ar[d]^{f_1} & A_2 \ar[r]^{a_2}\ar[d]^{f_2} & \cdots \ar[r]^{a_{n-1}}& A_{n} \ar[r]^{a_{n}\hspace{2mm}}\ar[d]^{f_{n}} & A_{n+1} \ar[r]^{a_{n+1}}\ar[d]^{f_{n+1}} & \Sigma A_0 \ar[d]^{\Sigma f_0}\\
B_0 \ar[r]^{b_0} & B_1 \ar[r]^{b_1} & B_2 \ar[r]^{b_2} & \cdots \ar[r]^{b_{n-1}}& B_{n} \ar[r]^{b_{n}\hspace{2mm}} & B_{n+1} \ar[r]^{b_{n+1}} & \Sigma B_0
}$$
commutes. It is an {\em isomorphism} if $f_0, f_1, \cdots, f_{n+1}$ are all isomorphisms in $\C$.

%We review the definition of $(n+2)$-angulated category, where it is mentioned in \cite{bt, ahbt} that both (N4) and (N4$^*$) can be the fourth axiom of $(n+2)$-angulated category, but no detailed proof is given, and Theorem \ref{th N4 equivalent} in this paper gives their equivalent proof.

Recently, Geiss, Keller, and Oppermann \cite{GKO} have discussed the axioms for an $(n+2)$-angulated category. In particular, they introduced a higher ``octahedral axiom" and showed that it is equivalent to the mapping cone axiom for an $(n+2)$-angulated category, as seen in \cite[Theorem 4.4]{bt}. Let's now recall the definition of an $(n+2)$-angulated category.

\begin{definition}\label{d1}\cite{GKO,bt}
An $(n+2)$-$angulated$ $category$ is a triple $(\C, \Sigma, \Phi)$, where $\C$ is an additive category, $\Sigma$ is an automorphism of $\C$, and $\Phi$ is a class of $(n+2)$-$\Sigma$-sequences (whose elements are called $(n+2)$-$angles$), which satisfies the following axioms:
\begin{itemize}[leftmargin=3em]
\item[\bf (N1)]
\begin{itemize}
\item[(a)] The class $\Phi$ is closed under isomorphisms, direct sums and  direct summands.

\item[(b)] For any object $A\in\C$, the following trivial sequence
$$A\xrightarrow{1} A\rightarrow 0\rightarrow 0\rightarrow\cdots\rightarrow 0\rightarrow \Sigma A$$
 belongs to $\Phi$.

\item[(c)]
Each morphism  $a_0\colon A_0\rightarrow A_1$ in $\C$, there exists an $(n+2)$-$\Sigma$-sequence in
$\Phi$ whose first morphism is $a_0$.

\end{itemize}
\item[\bf (N2)] An $(n+2)$-$\Sigma$-sequence belongs to $\Phi$ if and only if its left rotation belongs to $\Phi$.

\item[\bf (N3)] Given the solid part of the commutative diagram
$$\xymatrix{
A_0 \ar[r]^{a_0}\ar[d]^{f_0} & A_1 \ar[r]^{a_1}\ar[d]^{f_1} & A_2 \ar[r]^{a_2}\ar@{-->}[d]^{f_2} & \cdots \ar[r]^{a_{n-1}}& A_{n} \ar[r]^{a_{n}\hspace{2mm}}\ar@{-->}[d]^{f_{n}} & A_{n+1} \ar[r]^{a_{n+1}}\ar@{-->}[d]^{f_{n+1}} & \Sigma A_0 \ar[d]^{\Sigma f_0}\\
B_0 \ar[r]^{b_0} & B_1 \ar[r]^{b_1} & B_2 \ar[r]^{b_2} & \cdots \ar[r]^{b_{n-1}}& B_{n} \ar[r]^{b_{n}\hspace{2mm}} & B_{n+1} \ar[r]^{b_{n+1}} & \Sigma B_0
}$$
with rows in $\Phi$. Then there exist the dotted morphisms such that the above diagram commutes i.e. $(f_0, f_1, \cdots, f_{n+1})$ is a morphism of  $(n+2)$-$\Sigma$-sequences.

\item[\bf (N4)] Given a commutative diagram
$$\xymatrix{A_0 \ar[r]^{a_0}\ar[d]^{f_0} & A_1 \ar[r]^{a_1}\ar[d]^{f_1} & A_2 \ar[r]^{a_2} & \cdots \ar[r]^{a_{n-1}}& A_{n} \ar[r]^{a_{n}\hspace{2mm}} & A_{n+1} \ar[r]^{a_{n+1}} & \Sigma A_0 \ar[d]^{\Sigma f_0}\\
A_0 \ar[r]^{b_0} & B_1 \ar[r]^{b_1}\ar[d]^{c_1} & B_2 \ar[r]^{b_2} & \cdots \ar[r]^{b_{n-1}}& B_{n} \ar[r]^{b_{n}\hspace{2mm}} & B_{n+1} \ar[r]^{b_{n+1}} & \Sigma A_0\\
 & C_2 \ar[d]^{c_{2}}\\
 & \vdots \ar[d]^{c_{n-1}}\\
 & C_n\ar[d]^{c_n}\\
 & C_{n+1}\ar[d]^{c_{n+1}}\\
 & \Sigma A_1
 }$$
with the first two rows and the second column in $\Phi$. Then there exist morphisms
$$f_i: A_i \rightarrow B_i~~(i=2, 3, \cdots, n+1)$$

\newpage
$$\begin{aligned}
%&f_i: A_i \rightarrow B_i~~(i=2, 3, \cdots, n+1)\\
&g_i: B_i \rightarrow C_i~~(i=2, 3, \cdots, n+1)\\
&h_i: A_i \rightarrow C_{i-1}~~(i=3, 4, \cdots, n+1)
\end{aligned}$$
with the following two properties:
\begin{itemize}
\item[(a)] The sequence $(f_0, f_1, \cdots, f_{n+1})$ is a morphism of  $(n+2)$-$\Sigma$-sequences;

\item[(b)] The following $(n+2)$-$\Sigma$-sequence
$$A_2\xrightarrow{\left[
                    \begin{smallmatrix}
                      a_2 \\
                      f_2 \\
                    \end{smallmatrix}
                  \right]} A_3\oplus B_2\xrightarrow{\left[
                             \begin{smallmatrix}
                               -a_3 & 0 \\
                               f_3 & -b_2 \\
                               h_3 & g_2 \\
                             \end{smallmatrix}
                           \right]}
 A_4\oplus B_3\oplus C_2\xrightarrow{~\alpha_1~}A_5\oplus B_4\oplus C_3\xrightarrow{~\alpha_2~}\cdots\hspace{20mm}$$
$$\cdots\xrightarrow{~\alpha_{n-3}~}A_{n+1}\oplus B_{n}\oplus C_{n-1}\xrightarrow{~\beta~}B_{n+1}\oplus C_{n}\xrightarrow{\left[g_{n+1}~ c_{n}\right]}C_{n+1}\xrightarrow{\Sigma a_1\circ c_{n+1}}\Sigma A_2 $$
belongs to $\Phi$, where
$$\alpha_i=\left[\begin{matrix}
   -a_{i+3} & 0 & 0\\[1mm]
  (-1)^if_{i+3} & -b_{i+2} & 0\\[1mm]
   h_{i+3} & g_{i+2} & c_{i+1}\\
 \end{matrix}\right],~\beta=\left[\begin{matrix}
  (-1)^nf_{n+1} & -b_{n} & 0\\[1mm]
   h_{n+1} & g_{n} & c_{n-1}\\
 \end{matrix}\right],$$ and $c_{n+1}g_{n+1}=\Sigma a_0\circ b_{n+1}$.
\end{itemize}

\end{itemize}
\vspace{2mm}

Arentz-Hansen, Bergh and Thaule  mentioned the following in their paper \cite{ahbt}: ``The following axiom (N4$^\ast$) is not strictly the same as axiom (N4) in \cite{bt}. However, it follows from the proofs in \cite[Section 4]{bt} that the two are equivalent."
Unfortunately, there is no detailed proof provided, so for the convenience of the readers, we present a detailed proof.
\begin{itemize}[leftmargin=3.5em]
\item[\bf (N4$^{*}$)]  Given the solid part of the following commutative diagram
  $$
 \begin{tikzpicture}
       \draw[red] (-3.67,-0.68) node{$(\spadesuit)$};

      \diagram{d}{2.5em}{2.5em}{
        A_0 & A_1 & A_2 & A_3 & \cdots & A_{n} & A_{n+1} & \Sigma A_0\\
        A_0 & B_1 & B_2 & B_3 & \cdots & B_{n} & B_{n+1} & \Sigma
        A_0\\
        A_1 & B_1 & C_2 & C_3 & \cdots & C_{n} & C_{n+1} & \Sigma
        A_1\\
      };

      \path[->,midway,font=\scriptsize]
        (d-1-1) edge node[above] {$a_0$} (d-1-2)
        ([xshift=-0.1em] d-1-1.south) edge[-] ([xshift=-0.1em] d-2-1.north)
        ([xshift=0.1em] d-1-1.south) edge[-] ([xshift=0.1em] d-2-1.north)
        (d-1-2) edge node[above] {$a_1$} (d-1-3)
                     edge node[right] {$f_1$} (d-2-2)
        (d-1-3) edge node[above] {$a_2$} (d-1-4)
                     edge[densely dashed] node[right] {$f_2$} (d-2-3)
        (d-1-4) edge node[above] {$a_3$} (d-1-5)
                     edge[densely dashed] node[right] {$f_3$} (d-2-4)
        (d-1-5) edge node[above] {$a_{n-1}$} (d-1-6)
        (d-1-6) edge node[above] {$a_{n}$} (d-1-7)
                     edge[densely dashed] node[right] {$f_{n}$} (d-2-6)
        (d-1-7) edge node[above] {$a_{n+1}$} (d-1-8)
                     edge[densely dashed] node[right] {$f_{n+1}$} (d-2-7)
        ([xshift=-0.1em] d-1-8.south) edge[-] ([xshift=-0.1em] d-2-8.north)
        ([xshift=0.1em] d-1-8.south) edge[-] ([xshift=0.1em] d-2-8.north)
        (d-2-1) edge node[above] {$b_0$} (d-2-2)
                     edge node[right] {$a_0$} (d-3-1)
        (d-2-2) edge node[above] {$b_1$} (d-2-3)
        ([xshift=-0.1em] d-2-2.south) edge[-] ([xshift=-0.1em] d-3-2.north)
        ([xshift=0.1em] d-2-2.south) edge[-] ([xshift=0.1em] d-3-2.north)
        (d-2-3) edge node[above] {$b_2$} (d-2-4)
                     edge[densely dashed] node[right] {$g_2$} (d-3-3)
        (d-2-4) edge node[above] {$b_3$} (d-2-5)
                     edge[densely dashed] node[right] {$g_3$} (d-3-4)
        (d-2-5) edge node[above] {$b_{n - 1}$} (d-2-6)
        (d-2-6) edge node[above] {$b_{n}$} (d-2-7)
                     edge[densely dashed] node[right] {$g_{n}$} (d-3-6)
        (d-2-7) edge node[above] {$b_{n+1}$} (d-2-8)
                     edge[densely dashed] node[right] {$g_{n+1}$} (d-3-7)
        (d-2-8) edge node[right] {$\Sigma a_0$} (d-3-8)
        (d-3-1) edge node[above] {$f_1$} (d-3-2)
        (d-3-2) edge node[above] {$c_1$} (d-3-3)
        (d-3-3) edge node[above] {$c_2$} (d-3-4)
        (d-3-4) edge node[above] {$c_3$} (d-3-5)
        (d-3-5) edge node[above] {$c_{n - 1}$} (d-3-6)
        (d-3-6) edge node[above] {$c_{n}$} (d-3-7)
        (d-3-7) edge node[above] {$c_{n+1}$} (d-3-8)
        (d-1-4) edge[densely dashed,out=-102,in=30] node[pos=0.15,left] {$h_3$}
          (d-3-3)
        (d-1-7) edge[densely dashed,out=-102,in=30] node[pos=0.15,left] {$h_{n+1}$}
          (d-3-6);
    \end{tikzpicture}
$$

  with rows in $\Phi$, there exist the dotted
  morphisms such that each square commutes,  and the following
  $(n+2)$-$\Sigma$-sequence
$$A_2\xrightarrow{\left[
                    \begin{smallmatrix}
                      a_2 \\
                      f_2 \\
                    \end{smallmatrix}
                  \right]} A_3\oplus B_2\xrightarrow{\left[
                             \begin{smallmatrix}
                               -a_3 & 0 \\
                               f_3 & -b_2 \\
                               h_3 & g_2 \\
                             \end{smallmatrix}
                           \right]}
 A_4\oplus B_3\oplus C_2\xrightarrow{~\alpha_1~}A_5\oplus B_4\oplus C_3\xrightarrow{~\alpha_2~}\cdots\hspace{20mm}$$
$$\cdots\xrightarrow{~\alpha_{n-3}~}A_{n+1}\oplus B_{n}\oplus C_{n-1}\xrightarrow{~\beta~}B_{n+1}\oplus C_{n}\xrightarrow{\left[g_{n+1}~ c_{n}\right]}C_{n+1}\xrightarrow{\Sigma a_1\circ c_{n+1}}\Sigma A_2 $$
belongs to $\Phi$, where
$\alpha_i=\left[\begin{matrix}
   -a_{i+3} & 0 & 0\\[1mm]
  (-1)^if_{i+3} & -b_{i+2} & 0\\[1mm]
   h_{i+3} & g_{i+2} & c_{i+1}\\
 \end{matrix}\right],~\beta=\left[\begin{matrix}
  (-1)^nf_{n+1} & -b_{n} & 0\\[1mm]
   h_{n+1} & g_{n} & c_{n-1}\\
 \end{matrix}\right]$.
   \end{itemize}
\end{definition}

The following result show that axiom (N4) may be replaced by  axiom (N4$^\ast$).

\begin{theorem}\label{th N4 equivalent}
If $\Phi$ is a collection of $(n+2)$-$\Sigma$-sequences satisfying axioms \emph{(N1),~(N2)} and {\rm (N3)}, then $\Phi$ satisfies \emph{(N4)} if and only if
$\Phi$ satisfies \emph{(N4$^{*}$)}.
\end{theorem}

\proof {\bf Sufficiency.}  This is evident.

{\bf Necessity}. By (N4$^\ast$), we have the following
commutative diagram
\begin{equation}\label{(2.1) N4^* three rows diagram}
\begin{split}
\xymatrix{A_0 \ar[r]^{a_0} \ar@{=}[d] & A_1 \ar[r]^{a_1} \ar[d]^{f_1} & A_2 \ar[r]^{a_2} &\cdots \ar[r]^{a_{n-1}} & A_n \ar[r]^{a_n\hspace{2mm}} & A_{n+1} \ar[r]^{a_{n+1}} & \Sigma A_0 \ar@{=}[d]\\
A_0 \ar[r]^{b_0} \ar[d]^{a_0} & B_1 \ar[r]^{b_1} \ar@{=}[d] & B_2 \ar[r]^{b_2} &\cdots \ar[r]^{b_{n-1}} & B_n \ar[r]^{b_n\hspace{2mm}} & B_{n+1} \ar[r]^{b_{n+1}} & \Sigma A_0 \ar[d]^{\Sigma a_0}\\
A_1 \ar[r]^{f_1} & B_1 \ar[r]^{c_1} & C_2 \ar[r]^{c_2} &\cdots \ar[r]^{c_{n-1}} & C_n \ar[r]^{c_n\hspace{2mm}} & C_{n+1} \ar[r]^{c_{n+1}} & \Sigma A_1
}
\end{split}
\end{equation}
with rows in $\Phi$. By \cite[Lemma 4.1]{bt}, for the upper part of (\ref{(2.1) N4^* three rows diagram}), $f_2, f_3, \cdots, f_{n+1}$ can be chosen to complete the diagram of morphism
$$\xymatrix{A_0 \ar[r]^{a_0} \ar@{=}[d] & A_1 \ar[r]^{a_1} \ar[d]^{f_1} & A_2 \ar[r]^{a_2} \ar@{-->}[d]^{f_2} &\cdots \ar[r]^{a_{n-1}} & A_n \ar[r]^{a_n\hspace{2mm}} \ar@{-->}[d]^{f_n} & A_{n+1} \ar[r]^{a_{n+1}} \ar@{-->}[d]^{f_{n+1}} & \Sigma A_0 \ar@{=}[d]\\
A_0 \ar[r]^{b_0} & B_1 \ar[r]^{b_1} & B_2 \ar[r]^{b_2} &\cdots \ar[r]^{b_{n-1}} & B_n \ar[r]^{b_n\hspace{2mm}} & B_{n+1} \ar[r]^{b_{n+1}} & \Sigma A_0
}$$
such that the following $(n+2)$-$\Sigma$-sequence
$$A_1 \xrightarrow{\left[
                    \begin{smallmatrix}
                      -a_1 \\
                      f_1 \\
                    \end{smallmatrix}
                  \right]} A_2\oplus B_1\xrightarrow{\left[
                             \begin{smallmatrix}
                                a_2 & 0 \\
                                 f_2 & b_1 \\
                             \end{smallmatrix}
                           \right]} A_3\oplus B_2\xrightarrow{\left[
                             \begin{smallmatrix}
                                a_3 & 0 \\
                                 -f_3 & b_2 \\
                             \end{smallmatrix}
                           \right]} \cdots \hspace{50mm}$$
$$\cdots\xrightarrow{\left[
 \begin{smallmatrix}
     a_n & 0 \\
   (-1)^nf_n & b_{n-1} \\
  \end{smallmatrix}
  \right]} A_{n+1} \oplus B_n \xrightarrow{\left[
                            \begin{smallmatrix}
                                 (-1)^{n+1}f_{n+1} &~ b_{n} \\
                             \end{smallmatrix}
                             \right]} B_{n+1} \xrightarrow{\Sigma a_0\circ b_{n+1}}\Sigma A_1 $$
belongs to $\Phi$.
Note that we have the following solid commutative diagram
$$\xymatrix@C=1.5cm{A_1 \ar[r]^{\left[
                    \begin{smallmatrix}
                      -a_1 \\
                      f_1 \\
                    \end{smallmatrix}
                  \right]\hspace{6mm}} \ar@{=}[d] & A_2\oplus B_1 \ar[r]^{\left[
                             \begin{smallmatrix}
                                a_2 & 0 \\
                                 f_2 & b_1 \\
                             \end{smallmatrix}
                           \right]\hspace{2mm}} \ar[d]^{\left[
                             \begin{smallmatrix}
                                0 & 1 \\
                             \end{smallmatrix}
                           \right]} & A_3\oplus B_2 \ar[r]^{\hspace{3mm}\left[
                             \begin{smallmatrix}
                                a_3 & 0 \\
                                 -f_3 & b_2 \\
                             \end{smallmatrix}
                           \right]} \ar@{-->}[d]^{\left[
                             \begin{smallmatrix}
                                h_3 & g_2 \\
                             \end{smallmatrix}
                           \right]} & \cdots \\
A_1 \ar[r]^{f_1} & B_1 \ar[r]^{c_1} & C_2 \ar[r]^{c_2} & \cdots
}\hspace{40mm}$$
\begin{equation}\label{(2.2) mapping cone and C}
\begin{split}
\xymatrix@C=2.5cm{\cdots \ar[r]^{\left[
 \begin{smallmatrix}
     a_n & 0 \\
   (-1)^nf_n & b_{n-1} \\
  \end{smallmatrix}
  \right]\hspace{7mm}} & A_{n+1} \oplus B_n \ar[r]^{\hspace{4mm}\left[
                            \begin{smallmatrix}
                                 (-1)^{n+1}f_{n+1} & b_{n} \\
                             \end{smallmatrix}
                             \right]} \ar@{-->}[d]^{\left[
                             \begin{smallmatrix}
                                h_{n+1} & g_n \\
                             \end{smallmatrix}
                           \right]} & B_{n+1} \ar[r]^{\Sigma a_0 \circ b_{n+1}} \ar@{-->}[d]^{g_{n+1}} & \Sigma A_1 \ar@{=}[d]\\
\cdots \ar[r]^{c_{n-1}} & C_n \ar[r]^{c_n} & C_{n+1} \ar[r]^{c_{n+1}} & \Sigma A_1
}
\end{split}
\end{equation}
where rows in $\Phi$, by \cite[Lemma 4.1]{bt}, there exist morphisms $g_i: B_i \rightarrow C_i~~(i= 2, \cdots, n+1)$, and $h_i: A_i \rightarrow C_{i-1}~~(i= 3, \cdots, n+1)$ make the diagram (\ref{(2.2) mapping cone and C}) commutes, and the following mapping cone
$$A_2 \oplus B_1 \oplus A_1 \xrightarrow{\left[
                    \begin{smallmatrix}
                      -a_2 & 0 & 0\\
                      -f_2 & -b_1 & 0 \\
                      0 & 1 & f_1\\
                    \end{smallmatrix}
                  \right]} A_3 \oplus B_2\oplus B_1 \xrightarrow{~~\alpha_0~~} A_4\oplus B_3 \oplus C_2 \xrightarrow{~~\alpha_1~~} \cdots \hspace{50mm}$$
$$\cdots \xrightarrow{ \alpha_{n-3}} A_{n+1} \oplus B_n \oplus C_{n-1} \xrightarrow{\beta} B_{n+1} \oplus C_{n} \xrightarrow{\left[
                            \begin{smallmatrix}
                                -\Sigma a_0b_{n+1} & 0 \\
                                 g_{n+1} & c_{n} \\
                             \end{smallmatrix}
                             \right]}\Sigma A_1 \oplus C_{n+1} \xrightarrow{\left[
                            \begin{smallmatrix}
                                \Sigma a_1 & 0 \\
                                 -\Sigma f_1 & 0 \\
                                 1 & c_{n+1} \\
                             \end{smallmatrix}
                             \right]} \Sigma A_2 \oplus\Sigma B_1 \oplus \Sigma A_1$$
belongs to $\Phi$, where $\alpha_i=\left[\begin{matrix}
   -a_{i+3} & 0 & 0\\[1mm]
  (-1)^if_{i+3} & -b_{i+2} & 0\\[1mm]
   h_{i+3} & g_{i+2} & c_{i+1}\\
 \end{matrix}\right],~\beta=\left[\begin{matrix}
  (-1)^nf_{n+1} & -b_{n} & 0\\[1mm]
   h_{n+1} & g_{n} & c_{n-1}\\
 \end{matrix}\right].$

Since the diagram (\ref{(2.2) mapping cone and C}) commutes, we have the following
equalities.
$$
\begin{aligned}
& \left[
 \begin{matrix}
         h_3 & g_2 \\
              \end{matrix}
                \right]\left[
                    \begin{matrix}
                         a_2 & 0 \\
                        f_2 & b_1 \\
                          \end{matrix}
                  \right]=\left[
 \begin{matrix}
         h_3a_2+g_2f_2 & g_2b_1 \\
              \end{matrix}
                \right]=\left[
 \begin{matrix}
         0 & c_1 \\
              \end{matrix}
                \right]~~\Rightarrow~~ g_2b_1=c_1;\\
& \left[
 \begin{matrix}
         h_4 & g_3 \\
              \end{matrix}
                \right]\left[
                    \begin{matrix}
                         a_3 & 0 \\
                        -f_3 & b_2 \\
                          \end{matrix}
                  \right]=\left[
 \begin{matrix}
         h_4a_3+g_3f_3 & g_3b_2 \\
              \end{matrix}
                \right]=\left[
 \begin{matrix}
         c_2h_3 & c_2g_2 \\
              \end{matrix}
                \right]~~\Rightarrow~~ g_3b_2=c_2g_2;\\
& ~~~~~~~~~~~~~~~~~\cdots\\
& \left[
 \begin{matrix}
         (-1)^{n+1}g_{n+1}f_{n+1} & g_{n+1}b_n \\
              \end{matrix}
                \right]=\left[
 \begin{matrix}
         c_nh_{n+1} & c_ng_{n} \\
              \end{matrix}
                \right]~~\Rightarrow~~ g_{n+1}b_n=c_ng_{n};\\
&~c_{n+1}g_{n+1}=\Sigma a_0\circ b_{n+1}.
\end{aligned}$$
Thus we have the following commutative diagram
$$\xymatrix{A_0 \ar[r]^{b_0} \ar[d]^{a_0} & B_1 \ar[r]^{b_1} \ar@{=}[d] & B_2 \ar[r]^{b_2} \ar@{-->}[d]^{g_2} & \cdots \ar[r]^{b_{n-1}} & B_n \ar[r]^{b_n\hspace{2mm}} \ar@{-->}[d]^{g_n} & B_{n+1} \ar[r]^{b_{n+1}} \ar@{-->}[d]^{g_{n+1}} & \Sigma A_0 \ar[d]^{\Sigma a_0}\\
A_1 \ar[r]^{f_1} & B_1 \ar[r]^{c_1} & C_2 \ar[r]^{c_2} & \cdots \ar[r]^{c_{n-1}} & C_n \ar[r]^{c_{n}\hspace{2mm}} & C_{n+1} \ar[r]^{c_{n+1}} & \Sigma A_1.
}$$
Since $\Phi$ is closed under direct summands, then the following
commutative diagram
$$\xymatrix@C=2.5cm@R=1.5cm{A_2 \ar[r]^{\left[
                    \begin{smallmatrix}
                      a_2 \\
                      f_2 \\
                    \end{smallmatrix}
                  \right]} \ar[d]^{\left[
                    \begin{smallmatrix}
                      -1 \\
                      0 \\
                      0 \\
                    \end{smallmatrix}
                  \right]} & A_3\oplus B_2 \ar[r]^{\left[
                             \begin{smallmatrix}
                                -a_3 & 0 \\
                                 f_3 & -b_2 \\
                                 h_3 & g_2
                             \end{smallmatrix}
                           \right]\hspace{4mm}} \ar[d]^{\left[
                             \begin{smallmatrix}
                                1 & 0 \\
                                0 & 1 \\
                                0 & 0 \\
                             \end{smallmatrix}
                           \right]} & A_4\oplus B_3\oplus C_2 \ar@{=}[d] \\
A_2\oplus B_1\oplus A_1 \ar[r]^{\left[
                             \begin{smallmatrix}
                                -a_2 & 0 & 0 \\
                                 -f_2 & -b_1 & 0 \\
                                 0 & 1 & f_1 \\
                             \end{smallmatrix}
                           \right]} \ar[d]^{\left[
                             \begin{smallmatrix}
                                -1 & 0 & a_1 \\
                             \end{smallmatrix}
                           \right]} & A_3\oplus B_2\oplus B_1 \ar[r]^{\left[
                             \begin{smallmatrix}
                                -a_3 & 0 & 0 \\
                                 f_3 & -b_2 & 0 \\
                                 h_3 & g_2 & c_1 \\
                             \end{smallmatrix}
                           \right]} \ar[d]^{\left[
                             \begin{smallmatrix}
                                1 & 0 & 0 \\
                                0 & 1 & b_1\\
                             \end{smallmatrix}
                           \right]} & A_4\oplus B_3\oplus C_2 \ar@{=}[d]\\
A_2 \ar[r]^{\left[
                    \begin{smallmatrix}
                      a_2 \\
                      f_2 \\
                    \end{smallmatrix}
                  \right]} & A_3\oplus B_2 \ar[r]^{\left[
                             \begin{smallmatrix}
                                -a_3 & 0 \\
                                 f_3 & -b_2 \\
                                 h_3 & g_2
                             \end{smallmatrix}
                           \right]\hspace{4mm}} & A_4\oplus B_3\oplus C_2
}\hspace{25mm}$$
$$\xymatrix@C=3.5cm@R=1.5cm{\ar[r]^{\left[
                             \begin{smallmatrix}
                                -a_4 & 0 & 0 \\
                                 -f_4 & -b_3 & 0 \\
                                 h_4 & g_3 & c_2 \\
                             \end{smallmatrix}
                           \right]\hspace{3mm}} & \cdots \ar[r]^{\left[
 \begin{smallmatrix}
     (-1)^nf_{n+1} & -b_n & 0 \\
   h_{n+1} & g_n & c_{n-1} \\
  \end{smallmatrix}
  \right]\hspace{7mm}} & B_{n+1}\oplus C_n \ar@{=}[d]\\
\ar[r]^{\left[
                             \begin{smallmatrix}
                                -a_4 & 0 & 0 \\
                                 -f_4 & -b_3 & 0 \\
                                 h_4 & g_3 & c_2 \\
                             \end{smallmatrix}
                           \right]\hspace{3mm}} & \cdots \ar[r]^{\left[
                            \begin{smallmatrix}
                                 (-1)^nf_{n+1} & -b_n & 0 \\
                                 h_{n+1} & g_n & c_{n-1} \\
                             \end{smallmatrix}\right]\hspace{7mm}} & B_{n+1}\oplus C_n \ar@{=}[d]\\
\ar[r]^{\left[
                             \begin{smallmatrix}
                                -a_4 & 0 & 0 \\
                                 -f_4 & -b_3 & 0 \\
                                 h_4 & g_3 & c_2 \\
                             \end{smallmatrix}
                           \right]\hspace{3mm}} & \cdots
                              \ar[r]^{\left[
 \begin{smallmatrix}
     (-1)^nf_{n+1} & -b_n & 0 \\
   h_{n+1} & g_n & c_{n-1} \\
  \end{smallmatrix}
  \right]\hspace{7mm}} & B_{n+1}\oplus C_n
}$$
$$
\hspace{10mm}\xymatrix@C=2.5cm@R=1.5cm{\ar[r]^{\left[
                            \begin{smallmatrix}
                                 g_{n+1} & c_{n} \\
                             \end{smallmatrix}
                             \right]\hspace{4mm}} & C_{n+1} \ar[r]^{\Sigma a_1 \circ c_{n+1}} \ar[d]^{\left[
                            \begin{smallmatrix}
                                 -c_{n+1} \\
                                 1 \\
                             \end{smallmatrix}
                             \right]} & \Sigma A_2 \ar[d]^{\left[
                            \begin{smallmatrix}
                                 -1 \\
                                 0 \\
                                 0 \\
                             \end{smallmatrix}
                             \right]}\\
 \ar[r]^{\left[
                            \begin{smallmatrix}
                                 -\Sigma a_0\cdot b_{n+1} & 0\\
                                 g_{n+1} & c_n \\
                             \end{smallmatrix}
                             \right]\hspace{10mm}} & \Sigma A_1\oplus C_{n+1} \ar[r]^{\left[
                            \begin{smallmatrix}
                                 \Sigma a_1 & 0 \\
                                 -\Sigma f_1 & 0 \\
                                 1 & c_{n+1}\\
                             \end{smallmatrix}
                             \right]\hspace{6mm}} \ar[d]^{\left[
                            \begin{smallmatrix}
                                 0 & 1\\
                             \end{smallmatrix}
                             \right]} & \Sigma A_2\oplus \Sigma B_1\oplus \Sigma A_1 \ar[d]^{\left[
                            \begin{smallmatrix}
                                 -1 & 0 & \Sigma a_1\\
                             \end{smallmatrix}
                             \right]}\\
\ar[r]^{\left[
                            \begin{smallmatrix}
                                 g_{n+1} & c_{n} \\
                             \end{smallmatrix}
                             \right]\hspace{4mm}} & C_{n+1} \ar[r]^{\Sigma a_1\circ c_{n+1}} & \Sigma A_2
}
$$
shows that the following $(n+2)$-$\Sigma$-sequence
$$A_2 \xrightarrow{\left[
                    \begin{smallmatrix}
                      a_2 \\
                      f_2 \\
                    \end{smallmatrix}
                  \right]}A_3\oplus B_2 \xrightarrow{\left[
                             \begin{smallmatrix}
                                -a_3 & 0 \\
                                 f_3 & -b_2 \\
                                 h_3 & g_2
                             \end{smallmatrix}
                           \right]}A_4 \oplus B_3\oplus C_2\xrightarrow{\left[
                             \begin{smallmatrix}
                                -a_4 & 0 & 0 \\
                                 -f_4 & -b_3 & 0 \\
                                 h_4 & g_3 & c_2 \\
                             \end{smallmatrix}
                           \right]}\cdots\hspace{30mm}$$
$$\cdots \xrightarrow{\left[
 \begin{smallmatrix}
     (-1)^nf_{n+1} & -b_n & 0 \\
   h_{n+1} & g_n & c_{n-1} \\
  \end{smallmatrix}
  \right]}B_{n+1}\oplus C_n \xrightarrow{\left[
                            \begin{smallmatrix}
                                 g_{n+1} & c_n \\
                             \end{smallmatrix}
                             \right]}C_{n+1} \xrightarrow{\Sigma a_1 \circ c_{n+1}}\Sigma A_2$$
belongs to $\Phi$. This completes the proof.   \qed

\begin{remark}
In (N4), it's hard to see that the square $\red (\spadesuit)$ is commutative. With $\Sigma$ being an automorphism, it ensures that $\red (\spadesuit)$  is commutative; otherwise, it cannot be obtained.
\end{remark}

\section{Right $(n+2)$-angulated categories}
\setcounter{equation}{0}

Based on Theorem \ref{th N4 equivalent}, we refine the definition of a right $(n+2)$-angulated category, introducing some axioms that differ from the original definition \cite[Definition 2.1]{L2}.

\begin{definition}\label{right (n+2)-angulated category}
A \emph{right $(n+2)$-angulated category} is a triple $(\C, \Sigma, \Theta)$, where $\C$ is an additive category, $\Sigma$ is an endofunctor of $\C$, and $\Theta$ is a class of $(n+2)$-$\Sigma$-sequences (whose elements are called \emph{right $(n+2)$-angles}), which satisfies the following axioms:
\begin{itemize}[leftmargin=4em]
\item[\bf (RN1)]
\begin{itemize}
\item[(a)] The class $\Theta$ is closed under isomorphisms, direct sums and  direct summands.

\item[(b$^\ast$)] For any object $A\in\C$, the following trivial sequence
$$0\rightarrow A\xrightarrow{1} A\rightarrow 0\rightarrow \cdots\rightarrow 0$$
belongs to $\Theta$.

\item[(c)] Each morphism $a_0: A_0\rightarrow A_1$ in $\C$ can be embedded in a right $(n+2)$-angle, and $a_0$ is the first morphism.
\end{itemize}

\item[\bf (RN2)] If an $(n+2)$-$\Sigma$-sequence belongs to $\Theta$, then its left rotation belongs to $\Theta$.

\item[\bf (RN3)]{\bf (morphism axiom)} Given the solid part of the following commutative diagram
$$\xymatrix{
A_0 \ar[r]^{a_0}\ar[d]^{f_0} & A_1 \ar[r]^{a_1}\ar[d]^{f_1} & A_2 \ar[r]^{a_2}\ar@{-->}[d]^{f_2} & \cdots \ar[r]^{a_{n-1}}& A_{n} \ar[r]^{a_{n}}\ar@{-->}[d]^{f_{n}} & A_{n+1} \ar[r]^{a_{n+1}}\ar@{-->}[d]^{f_{n+1}} & \Sigma A_0 \ar[d]^{\Sigma f_0}\\
B_0 \ar[r]^{b_0} & B_1 \ar[r]^{b_1} & B_2 \ar[r]^{b_2} & \cdots \ar[r]^{b_{n-1}}& B_{n} \ar[r]^{b_{n}} & B_{n+1} \ar[r]^{b_{n+1}} & \Sigma B_0
}$$
with rows in $\Theta$. Then there exist the dotted morphisms such that the above diagram commutes. i.e. $(f_0, f_1, \cdots, f_{n+1})$ is a morphism of $(n+2)$-$\Sigma$-sequences.

\item[\bf (RN4$^\ast$)] {\bf (octahedral axiom)} Given the solid part of the following commutative diagram
$$
    \begin{tikzpicture}
      \diagram{d}{2.5em}{2.5em}{
        A_0 & A_1 & A_2 & A_3 & \cdots & A_{n} & A_{n+1} & \Sigma A_0\\
        A_0 & B_1 & B_2 & B_3 & \cdots & B_{n} & B_{n+1} & \Sigma
        A_0\\
        A_1 & B_1 & C_2 & C_3 & \cdots & C_{n} & C_{n+1} & \Sigma
        A_1\\
      };

      \path[->,midway,font=\scriptsize]
        (d-1-1) edge node[above] {$a_0$} (d-1-2)
        ([xshift=-0.1em] d-1-1.south) edge[-] ([xshift=-0.1em] d-2-1.north)
        ([xshift=0.1em] d-1-1.south) edge[-] ([xshift=0.1em] d-2-1.north)
        (d-1-2) edge node[above] {$a_1$} (d-1-3)
                     edge node[right] {$f_1$} (d-2-2)
        (d-1-3) edge node[above] {$a_2$} (d-1-4)
                     edge[densely dashed] node[right] {$f_2$} (d-2-3)
        (d-1-4) edge node[above] {$a_3$} (d-1-5)
                     edge[densely dashed] node[right] {$f_3$} (d-2-4)
        (d-1-5) edge node[above] {$a_{n-1}$} (d-1-6)
        (d-1-6) edge node[above] {$a_{n}$} (d-1-7)
                     edge[densely dashed] node[right] {$f_{n}$} (d-2-6)
        (d-1-7) edge node[above] {$a_{n+1}$} (d-1-8)
                     edge[densely dashed] node[right] {$f_{n+1}$} (d-2-7)
        ([xshift=-0.1em] d-1-8.south) edge[-] ([xshift=-0.1em] d-2-8.north)
        ([xshift=0.1em] d-1-8.south) edge[-] ([xshift=0.1em] d-2-8.north)
        (d-2-1) edge node[above] {$b_0$} (d-2-2)
                     edge node[right] {$a_0$} (d-3-1)
        (d-2-2) edge node[above] {$b_1$} (d-2-3)
        ([xshift=-0.1em] d-2-2.south) edge[-] ([xshift=-0.1em] d-3-2.north)
        ([xshift=0.1em] d-2-2.south) edge[-] ([xshift=0.1em] d-3-2.north)
        (d-2-3) edge node[above] {$b_2$} (d-2-4)
                     edge[densely dashed] node[right] {$g_2$} (d-3-3)
        (d-2-4) edge node[above] {$b_3$} (d-2-5)
                     edge[densely dashed] node[right] {$g_3$} (d-3-4)
        (d-2-5) edge node[above] {$b_{n - 1}$} (d-2-6)
        (d-2-6) edge node[above] {$b_{n}$} (d-2-7)
                     edge[densely dashed] node[right] {$g_{n}$} (d-3-6)
        (d-2-7) edge node[above] {$b_{n+1}$} (d-2-8)
                     edge[densely dashed] node[right] {$g_{n+1}$} (d-3-7)
        (d-2-8) edge node[right] {$\Sigma a_0$} (d-3-8)
        (d-3-1) edge node[above] {$f_1$} (d-3-2)
        (d-3-2) edge node[above] {$c_1$} (d-3-3)
        (d-3-3) edge node[above] {$c_2$} (d-3-4)
        (d-3-4) edge node[above] {$c_3$} (d-3-5)
        (d-3-5) edge node[above] {$c_{n - 1}$} (d-3-6)
        (d-3-6) edge node[above] {$c_{n}$} (d-3-7)
        (d-3-7) edge node[above] {$c_{n+1}$} (d-3-8)
        (d-1-4) edge[densely dashed,out=-102,in=30] node[pos=0.15,left] {$h_3$}
          (d-3-3)
        (d-1-7) edge[densely dashed,out=-102,in=30] node[pos=0.15,left] {$h_{n+1}$}
          (d-3-6);
    \end{tikzpicture}
$$
with commuting squares and rows in $\Theta$,  then there exist the dotted morphisms such that each square commutes,  and the following
$(n+2)$-$\Sigma$-sequence
$$A_2\xrightarrow{\left[
                    \begin{smallmatrix}
                      a_2 \\
                      f_2 \\
                    \end{smallmatrix}
                  \right]} A_3\oplus B_2\xrightarrow{\left[
                             \begin{smallmatrix}
                               -a_3 & 0 \\
                               f_3 & -b_2 \\
                               h_3 & g_2 \\
                             \end{smallmatrix}
                           \right]}
 A_4\oplus B_3\oplus C_2\xrightarrow{~\alpha_1~}A_5\oplus B_4\oplus C_3\xrightarrow{~\alpha_2~}\cdots\hspace{20mm}$$
$$\cdots\xrightarrow{~\alpha_{n-3}~}A_{n+1}\oplus B_{n}\oplus C_{n-1}\xrightarrow{~\beta~}B_{n+1}\oplus C_{n}\xrightarrow{\left[g_{n+1}~ c_{n}\right]}C_{n+1}\xrightarrow{\Sigma a_1\circ c_{n+1}}\Sigma A_2 $$
 belongs to $\Theta$,  where
$\alpha_i=\left[\begin{matrix}
   -a_{i+3} & 0 & 0\\[1mm]
  (-1)^if_{i+3} & -b_{i+2} & 0\\[1mm]
   h_{i+3} & g_{i+2} & c_{i+1}\\
 \end{matrix}\right],~\beta=\left[\begin{matrix}
  (-1)^nf_{n+1} & -b_{n} & 0\\[1mm]
   h_{n+1} & g_{n} & c_{n-1}\\
 \end{matrix}\right].$
  \end{itemize}
\end{definition}

\begin{remark}
Comparing \cite[Definition 2.1]{L2}, the conditions (RN1)(b$^\ast$) and (RN4$^\ast$) are different from those in \cite{L2}. Now we state those two conditions as presented in their paper. In \cite{L2}, (RN1)(b) takes the following form:
$A\xrightarrow{1} A\rightarrow 0\rightarrow 0\rightarrow\cdots\rightarrow 0\rightarrow \Sigma A$
and (RN4) is represented as: Given a commutative diagram
$$\xymatrix{A_0 \ar[r]^{a_0}\ar@{=}[d] & A_1 \ar[r]^{a_1}\ar[d]^{f_1} & A_2 \ar[r]^{a_2} & \cdots \ar[r]^{a_{n-1}}& A_{n} \ar[r]^{a_{n}\hspace{2mm}} & A_{n+1} \ar[r]^{a_{n+1}} & \Sigma A_0 \ar@{=}[d]\\
A_0 \ar[r]^{b_0} & B_1 \ar[r]^{b_1}\ar[d]^{c_1} & B_2 \ar[r]^{b_2} & \cdots \ar[r]^{b_{n-1}}& B_{n} \ar[r]^{b_{n}\hspace{2mm}} & B_{n+1} \ar[r]^{b_{n+1}} & \Sigma A_0\\
 & C_2 \ar[d]^{c_2}\\
 & \vdots \ar[d]^{c_{n-1}}\\
 & C_n\ar[d]^{c_n}\\
 & C_{n+1}\ar[d]^{c_{n+1}}\\
 & \Sigma A_1
 }$$
with the first two rows and the second column in $\Theta$. Then there exist morphisms
$$\begin{aligned}
&f_i: A_i \rightarrow B_i~~(i=2, 3, \cdots, n+1)\\
&g_i: B_i \rightarrow C_i~~(i=2, 3, \cdots, n+1)\\
&h_i: A_i \rightarrow C_{i-1}~~(i=3, 4, \cdots, n+1)
\end{aligned}$$
with the following two properties:
\begin{itemize}
\item[(a)] The sequence $(f_0, f_1, \cdots, f_{n+1})$ is a morphism of  $(n+2)$-$\Sigma$-sequences;

\item[(b)] The following $(n+2)$-$\Sigma$-sequence
$$A_2\xrightarrow{\left[
                    \begin{smallmatrix}
                      a_2 \\
                      f_2 \\
                    \end{smallmatrix}
                  \right]} A_3\oplus B_2\xrightarrow{\left[
                             \begin{smallmatrix}
                               -a_3 & 0 \\
                               f_3 & -b_2 \\
                               h_3 & g_2 \\
                             \end{smallmatrix}
                           \right]}
 A_4\oplus B_3\oplus C_2\xrightarrow{~\alpha_1~}A_5\oplus B_4\oplus C_3\xrightarrow{~\alpha_2~}\cdots\hspace{20mm}$$
$$\cdots\xrightarrow{~\alpha_{n-3}~}A_{n+1}\oplus B_{n}\oplus C_{n-1}\xrightarrow{~\beta~}B_{n+1}\oplus C_{n}\xrightarrow{\left[g_{n+1}~ c_{n}\right]}C_{n+1}\xrightarrow{\Sigma a_1\circ c_{n+1}}\Sigma A_2 $$
belongs to $\Theta$,
where
$\alpha_i=\left[\begin{matrix}
   -a_{i+3} & 0 & 0\\[1mm]
  (-1)^if_{i+3} & -b_{i+2} & 0\\[1mm]
   h_{i+3} & g_{i+2} & c_{i+1}\\
 \end{matrix}\right],~\beta=\left[\begin{matrix}
  (-1)^nf_{n+1} & -b_{n} & 0\\[1mm]
   h_{n+1} & g_{n} & c_{n-1}\\
 \end{matrix}\right]$,

  and $c_{n+1}g_{n+1}=\Sigma a_0\circ b_{n+1}$.
\end{itemize}
\end{remark}

\begin{remark}
(1) In \cite{BM}, the right triangulated category defines a trivial sequence of any object $A$ in $\C$ as follows:
$$0\rightarrow A\xrightarrow{1} A\rightarrow 0$$
Therefore, we define the trivial sequence in a right $(n+2)$-angulated category in the form given in (RN1)(b$^\ast$), and when (RN1)(b$^\ast$) holds, so does (RN1)(b).

(2) In Theorem \ref{th N4 equivalent},
 we know that if $\Phi$ satisfies (N1), (N2) and (N3), then (N4) and (N4$^*$) are equivalent. So we adopt  (N4$^*$)  to define the right $(n+2)$-angulated category.
 Note that (RN4) and (RN4$^\ast$) are not identical. It is easy to observe that (RN4$^\ast$) implies (RN4) is straightforward.
 Since $\Sigma$ is not an automorphism, it is challenging to ensure the reverse direction holds, making our definition stricter than the original one. Additionally, when $n=1$, our definition is in perfect alignment with the right triangulated category in the sense of
 Beligiannis and Marmaridis \cite[Definition 1.1]{BM}.

(3) If $\Sigma$ is an automorphism, and the condition (RN2) also holds the opposite.
In this case, the right $(n+2)$-angulated category is an $(n+2)$-angulated category in the sense of \cite{GKO,ahbt}.
\end{remark}

\section{Right $(n+2)$-angulated quotient categories}
In this section, we give some examples of $(n+2)$-angulated categories.
We need some preparations as follows.

Let $\C$ be an additive category and $f: A\rightarrow B$ is a morphism in $\C$.
Recall that $g:B\rightarrow C$ is a {\em weak cokernel} of $f$, if $gf=0$ and for each morphism $h:B\rightarrow X$ such that $hf=0$ there exists a morphism $s:C\rightarrow X$ (not necessarily unique) such that $sg=h$. These are reflected in the commutative diagram below
$$\xymatrix@C=1.5cm@R=1cm{A \ar[dr]_0 \ar[r]^{\hspace{4mm}f\quad\;\;}&B \ar[r]^{g\;\;}\ar[d]^{h}&C \ar@{-->}[dl]^{s}\\
      &X&}$$
Note that $g$ is a {\em cokernel} of $f$ if $g$ is an epimorphism. We can define {\em weak kernel} dually.

\begin{definition}\cite[Definition 2.2]{j}
Let $\C$ be an additive category and let $a_0: A_0\rightarrow A_1$ be a morphism in $\C$. An \emph{$n$-cokernel} of $a_0$ is a sequence
$$(a_1, \cdots, a_n):\xymatrix{A_1 \ar@{->}[r]^{a_1}&A_2 \ar@{->}[r]^{a_2} &\cdots \ar@{->}[r]^{a_{n-1}} & A_n \ar@{->}[r]^{a_n\hspace{2mm}}&A_{n+1}}$$
where the morphism $a_{k}$ is a weak cokernel of $a_{k-1}$ for all $1\leq k\leq n-1$ and $a_n$ is a cokernel of $a_{n-1}$. In this case, the sequence
$$\xymatrix{A_0 \ar[r]^{a_0} &A_1 \ar@{->}[r]^{a_1}&A_2 \ar@{->}[r]^{a_2} &\cdots \ar@{->}[r]^{a_{n-1}} & A_n \ar@{->}[r]^{a_n\hspace{2mm}}&A_{n+1}}$$
is called {\em right $n$-exact sequence}.
\end{definition}

\begin{definition}\cite[Definition 2.6]{L2}
Let $\X$ be a subcategory of $\C$ and let $a_0: A_0\rightarrow A_1$ be a morphism in $\C$. We say that $a_0$ has a \emph{special $n$-cokernel} with respect to $\X$, if $a_0$ has an \emph{$n$-cokernel}
$$(a_1,a_2, \cdots, a_n):\xymatrix{A_1 \ar@{->}[r]^{a_1}&X_2 \ar@{->}[r]^{a_2} &\cdots \ar@{->}[r]^{a_{n-1}} & X_n \ar@{->}[r]^{a_n\hspace{2mm}}&A_{n+1}}$$
where $X_2,X_3,\cdots,X_n\in \X$.
\end{definition}

\begin{lemma}\emph{\cite[Lemma 2.8]{L2}}\label{lemma 2.4 n-pushout diagram}
Let
$$\xymatrix{A_0 \ar[r]^{a_0} \ar@{=}[d] & A_1 \ar[r]^{a_1} \ar[d]^{f_1} & A_2 \ar[r]^{a_2} \ar[d]^{f_2} &\cdots \ar[r]^{a_{n-1}} & A_n \ar[r]^{a_n\hspace{2mm}} \ar[d]^{f_n} & A_{n+1} \ar[d]^{f_{n+1}}\\
A_0 \ar[r]^{b_0} & B_1 \ar[r]^{b_1} & B_2 \ar[r]^{b_2} &\cdots \ar[r]^{b_{n-1}} & B_n \ar[r]^{b_n\hspace{2mm}} & B_{n+1}
}$$
be a commutative diagram of right $n$-exact sequences. Then we have an $n$-pushout diagram as follows
$$\xymatrix{A_1 \ar[r]^{a_1} \ar[d]^{f_1} & A_2 \ar[r]^{a_2} \ar[d]^{f_2} &\cdots \ar[r]^{a_{n-1}} & A_n \ar[r]^{a_n\hspace{2mm}} \ar[d]^{f_n} & A_{n+1} \ar[d]^{f_{n+1}}\\
B_1 \ar[r]^{b_1} & B_2 \ar[r]^{b_2} &\cdots \ar[r]^{b_{n-1}} & B_n \ar[r]^{b_n\hspace{2mm}} & B_{n+1}.
}$$
That is to say,
$$A_1 \xrightarrow{\left[
                    \begin{smallmatrix}
                      -a_1 \\
                      f_1 \\
                    \end{smallmatrix}
                  \right]} A_2\oplus B_1\xrightarrow{\left[
                             \begin{smallmatrix}
                                -a_2 & 0 \\
                                 f_2 & b_1 \\
                             \end{smallmatrix}
                           \right]} A_3\oplus B_2\xrightarrow{\left[
                             \begin{smallmatrix}
                                -a_3 & 0 \\
                                 f_3 & b_2 \\
                             \end{smallmatrix}
                           \right]} \cdots \xrightarrow{\left[
 \begin{smallmatrix}
     -a_n & 0 \\
     f_n & b_{n-1} \\
  \end{smallmatrix}
  \right]} A_{n+1} \oplus B_n \xrightarrow{\left[
                            \begin{smallmatrix}
                                 f_{n+1} &~ b_{n} \\
                             \end{smallmatrix}
                             \right]} B_{n+1}$$
is a right $n$-exact sequence.
\end{lemma}

Let $\C$ be an additive category and $\X$ be a subcategory of $\C$. In this paper, when we say $\X$ is a subcategory of $\C$, it usually means that $\X$ is full, and is closed under isomorphisms, direct sums and direct summands.
In the quotient category $\C / \X$, the objects are the same as those in $\C$, and the morphisms are elements in $\C(A, B) / \X(A, B)$, where $\X(A, B)$ forms a subgroup of $\C(A, B)$, and these morphisms pass through objects within $\X$. We denote $\overline{f}$ as the image of the morphism $f: A\rightarrow B$ under the quotient functor $\C\rightarrow \C/\X$.

A morphism $f: A\rightarrow B$ in $\C$ is called $\X$-$monic$ if the sequence $\C(B, X) \xrightarrow{\C(f, X)} \C(A, X) \rightarrow 0$ is exact for any object $X$ in $\X$. Additionally, the morphism $f$ is referred to as a $left$ $\X$-$approximation$ of $A$ if $B$ belongs to $\X$.
The subcategory $\X$ is called \emph{covariantly finite} of $\C$, if any object $A \in \C$ has a left $\X$-approximation. We can define \emph{$\X$-epic, right $\X$-approximation and cintravariantly finite} dually.

\begin{definition}{\rm \cite[Proposition 3.1]{L2}}\label{proposition 2.5 functor}
Let $\C$ be an additive category, and $\X$ be a covariantly finite subcategory of $\C$. If every left $\X$-approximation has a special $n$-cokernel with respect to $\X$, then there exist an additive endofunctor $\Sigma: \C / \X \to \C / \X$, defined as follows:

For any object $A\in\C$, there exist a right $n$-exact sequence
$$A\xrightarrow{a_0}X_1\xrightarrow{a_1}X_2\xrightarrow{a_2}\cdots\xrightarrow{a_{n-1}} X_n\xrightarrow{a_n}B$$
where $a_0$ is a left $\X$-approximation of $A$ and $(a_1,a_2,\cdots, a_n)$ is a special $n$-cokernel of $a_0$.

For any morphism $f\colon A\rightarrow A'$, since $a_0$ is a left $\X$-approximation of $A$, we have the following commutative diagram
$$\xymatrix{
A\ar[r]^{a_0}\ar[d]^{f} & X_1\ar[r]^{a_1}\ar[d]^{f_1} & X_2\ar[r]^{a_2}\ar[d]^{f_2} & \cdots \ar[r]^{a_{n-1}} & X_{n}\ar[r]^{a_{n}}\ar[d]^{f_{n}}  & B \ar[d]^{g}\\
A'\ar[r]^{a_0'} & X'_1\ar[r]^{a_1'} & X'_2\ar[r]^{a_2'} & \cdots \ar[r]^{a_{n-1}'} & X'_{n}\ar[r]^{a_{n}'} &  B'.
}$$
Define an endofunctor $\Sigma\colon\C/\X\rightarrow\C/\X$ such that $\Sigma A=B$ and $\Sigma \overline{f}=\overline{g}$.
\end{definition}

\begin{definition}\cite[Definition 3.2]{L2}\label{definition 2.6 standard (n+2)-angle}
Let $\C$ be an additive category, and $\X$ be a covariantly finite subcategory of $\C$. Assume that every left $\X$-approximation has a special $n$-cokernel with respect to $\X$.
Let
$$\xymatrix{A_0 \ar[r]^{a_0} & A_1 \ar[r]^{a_1} & A_2 \ar[r]^{a_2} & \cdots \ar[r]^{a_{n-1}} & A_n \ar[r]^{a_n\hspace{2mm}} & A_{n+1}}$$
is a right $n$-exact sequence, where $a_0$ is $\X$-monic. Then there exist the following commutative diagram
$$\xymatrix{
A_0 \ar[r]^{a_0}\ar@{=}[d] & A_1 \ar[r]^{a_1}\ar[d]^{f_1} & A_2 \ar[r]^{a_2}\ar[d]^{f_2} & \cdots \ar[r]^{a_{n-1}}& A_{n} \ar[r]^{a_{n}\hspace{2mm}}\ar[d]^{f_{n}} & A_{n+1} \ar[d]^{f_{n+1}} \\
A_0 \ar[r]^{b_0} & X_1 \ar[r]^{b_1} & X_2 \ar[r]^{b_2} & \cdots \ar[r]^{b_{n-1}}& X_{n} \ar[r]^{b_{n}\hspace{2mm}} & \Sigma A_0
}$$
where $b_0$ is a left $\X$-approximation of $A$, $(b_1, b_2, \cdots, b_n)$ is a special $n$-cokernel of $b_0$. Then we have a complex
$$A_0 \xrightarrow{~\overline{a_0}~} A_1\xrightarrow{~\overline{a_1}~} A_2 \xrightarrow{~\overline{a_2}~} \cdots \xrightarrow{~\overline{a_n}~} A_{n+1} \xrightarrow{(-1)^n\overline{f_{n+1}}} \Sigma A_0.
$$
We define \emph{right $(n+2)$-angles} in $\C/\X$ as the complexes which are isomorphic to
complexes obtained in this way. We denote by $\Theta$ the class of right $(n+2)$-angles.
\end{definition}

\begin{lemma}\emph{\cite[Lemma 3.3]{L2}}\label{lemma 2.7 induced standard (n+2)-angle diagram}
Let
$$\xymatrix{A_0 \ar[r]^{a_0} \ar[d]^{f_0} & A_1 \ar[r]^{a_1} \ar[d]^{f_1} & A_2 \ar[r]^{a_2} \ar[d]^{f_2} &\cdots \ar[r]^{a_{n-1}} & A_n \ar[r]^{a_n\hspace{2mm}} \ar[d]^{f_n} & A_{n+1} \ar[d]^{f_{n+1}}\\
B_0 \ar[r]^{b_0} & B_1 \ar[r]^{b_1} & B_2 \ar[r]^{b_2} &\cdots \ar[r]^{b_{n-1}} & B_n \ar[r]^{b_n\hspace{2mm}} & B_{n+1}
}$$
be a commutative diagram of right $n$-exact sequences, where $a_0$ and $b_0$ are $\X$-monic. Then we have a commutative diagram
$$\xymatrix@C=1.5cm{A_0 \ar[r]^{\overline{a_0}} \ar[d]^{\overline{f_0}} & A_1 \ar[r]^{\overline{a_1}} \ar[d]^{\overline{f_1}} & A_2 \ar[r]^{\overline{a_2}} \ar[d]^{\overline{f_2}} &\cdots \ar[r]^{\overline{a_{n-1}}} & A_n \ar[r]^{\overline{a_n}\hspace{2mm}} \ar[d]^{\overline{f_n}} & A_{n+1} \ar[r]^{(-1)^n\overline{a_{n+1}}} \ar[d]^{\overline{f_{n+1}}} & \Sigma A_0 \ar[d]^{\Sigma\overline{f_0}}\\
B_0 \ar[r]^{\overline{b_0}} & B_1 \ar[r]^{\overline{b_1}} & B_2 \ar[r]^{\overline{b_2}} &\cdots \ar[r]^{\overline{b_{n-1}}} & B_n \ar[r]^{\overline{b_n}\hspace{2mm}} & B_{n+1} \ar[r]^{(-1)^n\overline{b_{n+1}}} & \Sigma B_0
}$$
of right $(n+2)$-angles in $\C/\X$.
\end{lemma}

The following result can be found in \cite[Theorem 3.4]{L2}.
Now we refine this result to satisfy our definition of a right $(n+2)$-angulated category
as stated in Definition \ref{right (n+2)-angulated category}.

\begin{theorem}\label{main2}
Let $\C$ be an additive category and $\X$ be a covariantly finite subcategory of $\C$. If every $\X$-monic morphism has an $n$-cokernel and every left $\X$-approximation has a special $n$-cokernel with respect to $\X$, then the quotient category $\C / \X$ is a right $(n+2)$-angulated category with respect to the endofunctor $\Sigma$ defined in Definition \emph{\ref{proposition 2.5 functor}} and right $(n+2)$-angles defined in Definition \emph{\ref{definition 2.6 standard (n+2)-angle}}.
\end{theorem}

\proof
We will now check the axioms of right $(n+2)$-angulated categories
as defined in Definition \ref{right (n+2)-angulated category}.
We only need to prove (RN1)(b$^*$) and (RN4$^*$), as the proofs for the other axioms are the same as the ones provided in \cite[Theorem 3.4]{L2} and can be omitted.

 The commutative diagram
$$\xymatrix{0 \ar[r] \ar@{=}[d] & A \ar[r]^{1} \ar[d] & A \ar[r] \ar[d] & 0 \ar[r] \ar[d] & \cdots \ar[r] & 0 \ar[r] \ar[d] &0 \ar@{=}[d]\\
0 \ar[r] & 0 \ar[r] & 0 \ar[r] &0 \ar[r] & \cdots \ar[r] & 0 \ar[r] & 0
}$$
shows that
$0\rightarrow A\xrightarrow{1} A\rightarrow 0\rightarrow \cdots\rightarrow 0$
belongs to $\Theta$. Thus (RN1)(b$^*$) holds.

Now we prove that (RN4$^*$) holds. Given the solid part of the diagram
$$
\xymatrix@C=1.5cm{A_0 \ar[r]^{\overline{a_0}} \ar@{=}[d] & A_1 \ar[r]^{\overline{a_1}} \ar[d]^{\overline{f_1}} & A_2 \ar[r]^{\overline{a_2}} &\cdots  \ar[r]^{\overline{a_n}} & A_{n+1} \ar[r]^{(-1)^n\overline{a_{n+1}}} & \Sigma A_0 \ar@{=}[d]\\
A_0 \ar[r]^{\overline{f_1a_0}} \ar[d]^{\overline{a_0}} & B_1 \ar[r]^{\overline{b_1}} \ar@{=}[d] & B_2 \ar[r]^{\overline{b_2}} &\cdots \ar[r]^{\overline{b_n}} & B_{n+1} \ar[r]^{(-1)^n\overline{b_{n+1}}} & \Sigma A_0 \ar[d]^{\Sigma \overline{a_0}}\\
A_1 \ar[r]^{\overline{f_1}} & B_1 \ar[r]^{\overline{c_1}} & C_2 \ar[r]^{\overline{c_2}} &\cdots \ar[r]^{\overline{c_n}} & C_{n+1} \ar[r]^{(-1)^n\overline{c_{n+1}}} & \Sigma A_1
}
$$
with commuting squares and with rows in $\Theta$.
Without loss of generality, we can assume that $a_0$ and $f_1$ are $\X$-monic.
 Thus we obtain that $f_1a_0$ is also $\X$-monic.
Moreover, we have the following solid part commutative diagram
\begin{equation}\label{right n-exact sequence}
\begin{split}
\xymatrix{
A_0 \ar[r]^{a_0}\ar@{=}[d] & A_1 \ar[r]^{a_1}\ar[d]^{f_1} & A_2 \ar[r]^{a_2}\ar@{-->}[d]^{f_2} & \cdots \ar[r]^{a_{n-1}}& A_{n} \ar[r]^{a_{n}\hspace{2mm}}\ar@{-->}[d]^{f_{n}} & A_{n+1} \ar@{-->}[d]^{f_{n+1}}\\
A_0 \ar[r]^{f_1a_0} & B_1 \ar[r]^{b_1} & B_2 \ar[r]^{b_2} & \cdots \ar[r]^{b_{n-1}}& B_{n} \ar[r]^{b_{n}\hspace{2mm}} & B_{n+1}
}
\end{split}
\end{equation}
where rows are right $n$-exact sequences. By the property of weak cokernels, there exist dotted morphisms such that the diagram (\ref{right n-exact sequence}) commutes.
By Lemma \ref{lemma 2.7 induced standard (n+2)-angle diagram}, we obtain the following commutative diagram
$$\xymatrix@C=1.5cm{A_0 \ar[r]^{\overline{a_0}} \ar@{=}[d] & A_1 \ar[r]^{\overline{a_1}} \ar[d]^{\overline{f_1}} & A_2 \ar[r]^{\overline{a_2}} \ar[d]^{\overline{f_2}} &\cdots \ar[r]^{\overline{a_{n-1}}} & A_n \ar[r]^{\overline{a_n}\hspace{2mm}} \ar[d]^{\overline{f_n}} & A_{n+1} \ar[r]^{(-1)^n\overline{a_{n+1}}} \ar[d]^{\overline{f_{n+1}}} & \Sigma A_0 \ar@{=}[d]\\
A_0 \ar[r]^{\overline{f_1a_0}} & B_1 \ar[r]^{\overline{b_1}} & B_2 \ar[r]^{\overline{b_2}} &\cdots \ar[r]^{\overline{b_{n-1}}} & B_n \ar[r]^{\overline{b_n}\hspace{2mm}} & B_{n+1} \ar[r]^{(-1)^n\overline{b_{n+1}}} & \Sigma A_0
}$$
 of right $(n+2)$-angles.
By the diagram (\ref{right n-exact sequence}) and Lemma \ref{lemma 2.4 n-pushout diagram}, we obtain the following right $n$-exact sequence
$$A_1\xrightarrow{\left[
                    \begin{smallmatrix}
                      -a_1 \\
                      f_1 \\
                    \end{smallmatrix}
                  \right]} A_2\oplus B_1\xrightarrow{\left[
                             \begin{smallmatrix}
                               -a_2 & 0 \\
                               f_2 & b_1 \\
                             \end{smallmatrix}
                           \right]}
 A_3\oplus B_2\xrightarrow{\left[
                             \begin{smallmatrix}
                               -a_3 & 0 \\
                               f_3 & b_2 \\
                             \end{smallmatrix}
                           \right]} \cdots \xrightarrow{\left[
                             \begin{smallmatrix}
                               -a_n & 0 \\
                               f_n & b_{n-1} \\
                             \end{smallmatrix}
                           \right]}A_{n+1}\oplus B_n \xrightarrow{\left[
                             \begin{smallmatrix}
                               f_{n+1} & b_{n} \\
                             \end{smallmatrix}
                           \right]}B_{n+1}.$$
The commutative diagram below{\small
$$
\xymatrix@C=1.7cm@R=1.5cm{A_1 \ar[r]^{\left[
                    \begin{smallmatrix}
                      -a_1 \\
                      f_1 \\
                    \end{smallmatrix}
                  \right]\hspace{6mm}} \ar@{=}[d] & A_2\oplus B_1 \ar[r]^{\left[
                             \begin{smallmatrix}
                               -a_2 & 0 \\
                               f_2 & b_1 \\
                             \end{smallmatrix}
                           \right]\hspace{1mm}} \ar@{=}[d] & A_3\oplus B_2 \ar[r]^{\hspace{3mm}\left[
                             \begin{smallmatrix}
                               -a_3 & 0 \\
                               f_3 & b_2 \\
                             \end{smallmatrix}
                           \right]} \ar[d]^{\left[
                             \begin{smallmatrix}
                               -1 & 0 \\
                               0 & 1 \\
                             \end{smallmatrix}
                           \right]} & \cdots\\
A_1 \ar[r]^{\left[
                    \begin{smallmatrix}
                      -a_1 \\
                      f_1 \\
                    \end{smallmatrix}
                  \right]\hspace{6mm}} & A_2\oplus B_1 \ar[r]^{\left[
                             \begin{smallmatrix}
                               a_2 & 0 \\
                               f_2 & b_1 \\
                             \end{smallmatrix}
                           \right]}& A_3\oplus B_2 \ar[r]^{\hspace{3mm}\left[
                             \begin{smallmatrix}
                               a_3 & 0 \\
                               -f_3 & b_2 \\
                             \end{smallmatrix}
                           \right]}& \cdots
}\hspace{45mm}$$
$$\hspace{45mm}\xymatrix@C=2.5cm@R=1.5cm{\cdots \ar[r]^{\left[
                             \begin{smallmatrix}
                               -a_n & 0 \\
                               f_n & b_{n-1} \\
                             \end{smallmatrix}
                           \right]\hspace{6mm}}& A_{n+1}\oplus B_n \ar[r]^{\hspace{3mm}\left[
                             \begin{smallmatrix}
                               f_{n+1} & b_{n} \\
                             \end{smallmatrix}
                           \right]} \ar[d]^{\left[
                             \begin{smallmatrix}
                               (-1)^{n+1} & 0 \\
                               0 & 1 \\
                             \end{smallmatrix}
                           \right]} & B_{n+1} \ar@{=}[d]\\
\cdots \ar[r]^{\left[
                             \begin{smallmatrix}
                               a_n & 0 \\
                               (-1)^nf_n & b_{n-1} \\
                             \end{smallmatrix}
                           \right]\hspace{7mm}}& A_{n+1}\oplus B_n \ar[r]^{\hspace{3mm}\left[
\begin{smallmatrix}
                               (-1)^{n+1}f_{n+1} & b_{n} \\
                             \end{smallmatrix}
                           \right]}& B_{n+1}
}$$
shows that
$$A_1\xrightarrow{\left[
                    \begin{smallmatrix}
                      -a_1 \\
                      f_1 \\
                    \end{smallmatrix}
                  \right]} A_2\oplus B_1\xrightarrow{\left[
                             \begin{smallmatrix}
                               a_2 & 0 \\
                               f_2 & b_1 \\
                             \end{smallmatrix}
                           \right]}
 A_3\oplus B_2\xrightarrow{\left[
                             \begin{smallmatrix}
                               a_3 & 0 \\
                               -f_3 & b_2 \\
                             \end{smallmatrix}
                           \right]} \cdots \xrightarrow{\left[
                             \begin{smallmatrix}
                               a_n & 0 \\
                               (-1)^nf_n & b_{n-1} \\
                             \end{smallmatrix}
                           \right]}A_{n+1}\oplus B_n \xrightarrow{\left[
                             \begin{smallmatrix}
                               (-1)^{n+1}f_{n+1} & b_{n} \\
                             \end{smallmatrix}
                           \right]}B_{n+1}$$}
is a right $n$-exact sequence. Hence we have the following solid part commutative diagram
{\small $$
\xymatrix@C=1.7cm{A_1 \ar[r]^{\left[
                    \begin{smallmatrix}
                      -a_1 \\
                      f_1 \\
                    \end{smallmatrix}
                  \right]\hspace{4mm}} \ar@{=}[d] & A_2\oplus B_1 \ar[r]^{\left[
                             \begin{smallmatrix}
                               a_2 & 0 \\
                               f_2 & b_1 \\
                             \end{smallmatrix}
                           \right]} \ar[d]^{\left[
                             \begin{smallmatrix}
                               0 & 1 \\
                             \end{smallmatrix}
                           \right]} & A_3\oplus B_2 \ar[r]^{\hspace{3mm}\left[
                             \begin{smallmatrix}
                               a_3 & 0 \\
                               -f_3 & b_2 \\
                             \end{smallmatrix}
                           \right]} \ar@{-->}[d]^{\left[
                             \begin{smallmatrix}
                               h_3 & g_2 \\
                             \end{smallmatrix}
                           \right]} & \cdots\\
A_1 \ar[r]^{f_1} & B_1 \ar[r]^{c_1} & C_2 \ar[r]^{c_2} & \cdots
}\hspace{45mm}$$}
\begin{equation}\label{mapping cone induces a direct and C}
\begin{split}
\hspace{45mm}\xymatrix@C=2.5cm{\cdots \ar[r]^{\left[
                             \begin{smallmatrix}
                               a_n & 0 \\
                               (-1)^nf_n & b_{n-1} \\
                             \end{smallmatrix}
                           \right]\hspace{7mm}}& A_{n+1}\oplus B_n \ar[r]^{\hspace{3mm}\left[
                             \begin{smallmatrix}
                               (-1)^{n+1}f_{n+1} & b_{n} \\
                             \end{smallmatrix}
                           \right]} \ar@{-->}[d]^{\left[
                             \begin{smallmatrix}
                               h_{n+1} & g_n \\
                             \end{smallmatrix}
                           \right]} & B_{n+1} \ar@{-->}[d]^{g_{n+1}}\\
\cdots \ar[r]^{c_{n-1}}& C_n \ar[r]^{c_n}& C_{n+1}
}
\end{split}
\end{equation}
with rows are right $n$-exact sequences. By the property of weak cokernels, there exist $g_i: B_i\rightarrow C_i~~(i=2, \cdots, n+1), h_i: A_i\rightarrow C_{i-1}~~(i=3, \cdots, n+1)$ such that the diagram (\ref{mapping cone induces a direct and C}) commutes.
Since the diagram (\ref{mapping cone induces a direct and C}) is commutative,
we have the follow equalities:
$$
\begin{aligned}
& \left[
 \begin{matrix}
         h_3 & g_2 \\
              \end{matrix}
                \right]\left[
                    \begin{matrix}
                         a_2 & 0 \\
                        f_2 & b_1 \\
                          \end{matrix}
                  \right]=\left[
 \begin{matrix}
         h_3a_2+g_2f_2 & g_2b_1 \\
              \end{matrix}
                \right]=\left[
 \begin{matrix}
         0 & c_1 \\
              \end{matrix}
                \right]~~\Rightarrow~~ g_2b_1=c_1;\\
& \left[
 \begin{matrix}
         h_4 & g_3 \\
              \end{matrix}
                \right]\left[
                    \begin{matrix}
                         a_3 & 0 \\
                        -f_3 & b_2 \\
                          \end{matrix}
                  \right]=\left[
 \begin{matrix}
         h_4a_3-g_3f_3 & g_3b_2 \\
              \end{matrix}
                \right]=\left[
 \begin{matrix}
         c_2h_3 & c_2g_2 \\
              \end{matrix}
                \right]~~\Rightarrow~~ g_3b_2=c_2g_2;\\
& ~~~~~~~~~~~~~~~~~\cdots\\
& \left[
 \begin{matrix}
         (-1)^{n+1}g_{n+1}f_{n+1} & g_{n+1}b_n \\
              \end{matrix}
                \right]=\left[
 \begin{matrix}
         c_nh_{n+1} & c_ng_{n} \\
              \end{matrix}
                \right]~~\Rightarrow~~ g_{n+1}b_n=c_ng_{n}.
\end{aligned}$$
Thus we obtain the following commutative diagram
$$
\xymatrix{
A_0 \ar[r]^{f_1a_0}\ar[d]^{a_0} & B_1 \ar[r]^{b_1}\ar@{=}[d] & B_2 \ar[r]^{b_2}\ar[d]^{g_2} & \cdots \ar[r]^{b_{n-1}}& B_{n} \ar[r]^{b_{n}\hspace{2mm}}\ar[d]^{g_{n}} & B_{n+1} \ar[d]^{g_{n+1}}\\
A_1 \ar[r]^{f_1} & B_1 \ar[r]^{c_1} & C_2 \ar[r]^{c_2} & \cdots \ar[r]^{c_{n-1}}& C_{n} \ar[r]^{c_{n}\hspace{2mm}} & C_{n+1}.
}
$$
By Lemma \ref{lemma 2.7 induced standard (n+2)-angle diagram}, we get the following commutative diagram
$$\xymatrix@C=1.5cm{A_0 \ar[r]^{\overline{f_1a_0}} \ar[d]^{\overline{a_0}} & B_1 \ar[r]^{\overline{b_1}} \ar@{=}[d] & B_2 \ar[r]^{\overline{b_2}} \ar@{->}[d]^{\overline{g_2}} &\cdots \ar[r]^{\overline{b_{n-1}}} & B_n \ar[r]^{\overline{b_n}} \ar@{->}[d]^{\overline{g_n}} & B_{n+1} \ar[r]^{(-1)^n\overline{b_{n+1}}} \ar@{->}[d]^{\overline{g_{n+1}}} & \Sigma A_0 \ar[d]^{\Sigma \overline{a_0}}\\
A_1 \ar[r]^{\overline{f_1}} & B_1 \ar[r]^{\overline{c_1}} & C_2 \ar[r]^{\overline{c_2}} &\cdots \ar[r]^{\overline{c_{n-1}}} & C_n \ar[r]^{\overline{c_n}} & C_{n+1} \ar[r]^{(-1)^n\overline{c_{n+1}}} & \Sigma A_1
}$$
of right $(n+2)$-angles.
By Lemma \ref{lemma 2.4 n-pushout diagram} for the diagram (\ref{mapping cone induces a direct and C}), we obtain a right $n$-exact sequence as follows
$$
A_2 \oplus B_1\xrightarrow{\left[
                    \begin{smallmatrix}
                      -a_2 & 0 \\
                      -f_2 & -b_1\\
                      0 & 1 \\
                    \end{smallmatrix}
                  \right]} A_3\oplus B_2\oplus B_1\xrightarrow{~\alpha_0~}
 A_4\oplus B_3\oplus C_2\xrightarrow{~\alpha_1~}A_5\oplus B_4\oplus C_3\xrightarrow{~\alpha_2~}\cdots\hspace{20mm}$$
\begin{equation}\label{A2B1 mapping cone}
\begin{split}
\hspace{20mm}\cdots \xrightarrow{~\alpha_{n-3}~}A_{n+1}\oplus B_{n}\oplus C_{n-1}\xrightarrow{~\beta~}B_{n+1}\oplus C_{n}\xrightarrow{\left[g_{n+1}~ c_{n}\right]}C_{n+1}
\end{split}
\end{equation}
where
$\alpha_i=\left[\begin{matrix}
   -a_{i+3} & 0 & 0\\[1mm]
  (-1)^if_{i+3} & -b_{i+2} & 0\\[1mm]
   h_{i+3} & g_{i+2} & c_{i+1}\\
 \end{matrix}\right],~~\beta=\left[\begin{matrix}
  (-1)^nf_{n+1} & -b_{n} & 0\\[1mm]
   h_{n+1} & g_{n} & c_{n-1}\\
 \end{matrix}\right]$.

Next we want to show that the following sequence is a right $n$-exact sequence
\begin{equation}\label{RN4 last direct sums diagram induce}
A_2 \xrightarrow{\left[
                    \begin{smallmatrix}
                      a_2 \\
                      f_2 \\
                    \end{smallmatrix}
                  \right]}A_3\oplus B_2 \xrightarrow{\left[
                             \begin{smallmatrix}
                                -a_3 & 0 \\
                                 f_3 & -b_2 \\
                                 h_3 & g_2
                             \end{smallmatrix}
                           \right]}A_4 \oplus B_3\oplus C_2\xrightarrow{\left[
                             \begin{smallmatrix}
                                -a_4 & 0 & 0 \\
                                 -f_4 & -b_3 & 0 \\
                                 h_4 & g_3 & c_2 \\
                             \end{smallmatrix}
                           \right]}\cdots\hspace{40mm}
\end{equation}
$$
\cdots \xrightarrow{\left[
 \begin{smallmatrix}
     -a_n & 0 & 0 \\
   (-1)^{n-1}f_n & -b_{n-1} & 0 \\
   h_n & g_{n-1} & c_{n-2}\\
  \end{smallmatrix}
  \right]}A_{n+1}\oplus B_n\oplus C_{n-1} \xrightarrow{\left[
                            \begin{smallmatrix}
                                 (-1)^nf_{n+1} & -b_n & 0 \\
                                 h_{n+1} & g_n & c_{n-1}\\
                             \end{smallmatrix}
                             \right]}B_{n+1}\oplus C_n \xrightarrow{\left[
                            \begin{smallmatrix}
                                 g_{n+1} & c_n \\
                             \end{smallmatrix}
                             \right]}C_{n+1}.
$$
It suffices to prove that $\left[
                             \begin{matrix}
                                -a_3 & 0 \\
                                 f_3 & -b_2 \\
                                 h_3 & g_2
                             \end{matrix}
                           \right]$ is a weak cokernel of $\left[
                    \begin{matrix}
                      a_2 \\
                      f_2 \\
                    \end{matrix}
                  \right]$ and $\left[
                             \begin{matrix}
                                -a_4 & 0 & 0 \\
                                 -f_4 & -b_3 & 0 \\
                                 h_4 & g_3 & c_2 \\
                             \end{matrix}
                           \right]$ is a weak cokernel of $\left[
                             \begin{matrix}
                                -a_3 & 0 \\
                                 f_3 & -b_2 \\
                                 h_3 & g_2
                             \end{matrix}
                           \right]$.

Let $\left[
                    \begin{matrix}
                      s & t \\
                    \end{matrix}
                  \right]: A_3\oplus B_2\rightarrow M$ such that $\left[
                    \begin{matrix}
                      s & t \\
                    \end{matrix}
                  \right]\left[
                    \begin{matrix}
                      a_2 \\
                      f_2 \\
                    \end{matrix}
                  \right]=0$, i.e. $sa_2+tf_2=0$.
Then we have $\left[
\begin{matrix}
   s & t & tb_1\\
      \end{matrix}
       \right]: A_3\oplus B_2\oplus B_1\rightarrow M$ such that $\left[
\begin{matrix}
   s & t & tb_1\\
      \end{matrix}
       \right]\left[
        \begin{matrix}
          -a_2 & 0 \\
             -f_2 & -b_1\\
                0 & 1 \\
                 \end{matrix}
\right]=0$, we know that the sequence (\ref{A2B1 mapping cone}) is a right $n$-exact sequence, there exist a morphism $\left[
  \begin{matrix}
     p & q & r\\
      \end{matrix}
\right]: A_4\oplus B_3\oplus C_2\rightarrow M$ such that
$$\left[
  \begin{matrix}
     p & q & r\\
      \end{matrix}
\right]\left[
                             \begin{matrix}
                                -a_3 & 0 & 0\\
                                 f_3 & -b_2 & 0\\
                                 h_3 & g_2 & c_1
                             \end{matrix}
                           \right]=\left[
\begin{matrix}
   s & t & tb_1\\
      \end{matrix}
 \right].$$
 It follows that $-pa_3+qf_3+rh_3=s, -qb_2+rg_2=t, rc_1=tb_1$.
Then we have $$\left[
  \begin{matrix}
     p & q & r\\
      \end{matrix}
\right]\left[
                             \begin{matrix}
                                -a_3 & 0 \\
                                 f_3 & -b_2 \\
                                 h_3 & g_2
                             \end{matrix}
                           \right]=\left[
\begin{matrix}
   s & t \\
      \end{matrix}
 \right].$$
This is reflected in the commutative diagram below
$$\xymatrix@C=2cm@R=1.2cm{A_2\ar[dr]_0\ar[r]^{\left[
                    \begin{smallmatrix}
                      a_2 \\
                      f_2 \\
                    \end{smallmatrix}
                  \right]\quad\;\;}&A_3\oplus B_2\ar[r]^{\left[
                             \begin{smallmatrix}
                                -a_3 & 0 \\
                                 f_3 & -b_2 \\
                                 h_3 & g_2
                             \end{smallmatrix}
                           \right]\;\;\hspace{3mm}}\ar[d]^{\left[
\begin{smallmatrix}
   s & t \\
      \end{smallmatrix}
 \right]}&A_4\oplus B_3\oplus C_2\ar@{.>}[dl]^{\hspace{5mm}\left[
  \begin{smallmatrix}
     p & q & r\\
      \end{smallmatrix}
\right]}\\
      &M&}$$
This shows that $\left[
                             \begin{matrix}
                                -a_3 & 0 \\
                                 f_3 & -b_2 \\
                                 h_3 & g_2
                             \end{matrix}
                           \right]$ is a weak cokernel of $\left[
                    \begin{matrix}
                      a_2 \\
                      f_2 \\
                    \end{matrix}
                  \right]$.

Let $\left[
                    \begin{matrix}
                      u & v & w \\
                    \end{matrix}
                  \right]: A_4\oplus B_3\oplus C_2\rightarrow X$ such that $\left[
                    \begin{matrix}
                      u & v & w \\
                    \end{matrix}
                  \right]\left[
                             \begin{matrix}
                                -a_3 & 0 \\
                                 f_3 & -b_2 \\
                                 h_3 & g_2
                             \end{matrix}
                           \right]=0$, i.e. $-ua_3+vf_3+wh_3=0, -vb_2+wg_2=0$.
Then $-vb_2b_1+wg_2b_1=0$, i.e. $wc_1=0$. Thus we have $$\left[
                    \begin{matrix}
                      u & v & w \\
                    \end{matrix}
                  \right]\left[
                             \begin{matrix}
                                -a_3 & 0 & 0\\
                                 f_3 & -b_2 & 0\\
                                 h_3 & g_2 & c_1
                             \end{matrix}
                           \right]=0,$$
we know that the sequence (\ref{A2B1 mapping cone}) is a right $n$-exact sequence, there exist a morphism $\left[
  \begin{matrix}
     l & m & n\\
      \end{matrix}
\right]: A_5\oplus B_4\oplus C_3\rightarrow X$ such that $\left[
  \begin{matrix}
     l & m & n\\
      \end{matrix}
\right]\left[
                             \begin{matrix}
                                -a_4 & 0 & 0\\
                                 -f_4 & -b_3 & 0\\
                                 h_4 & g_3 & c_2
                             \end{matrix}
                           \right]=\left[
\begin{matrix}
   u & v & w\\
      \end{matrix}
 \right]$ i.e. we have the following commutative diagram
$$\xymatrix@C=2.5cm@R=1.2cm{A_3\oplus B_2\ar[dr]_0\ar[r]^{\left[
                             \begin{smallmatrix}
                                -a_3 & 0 \\
                                 f_3 & -b_2 \\
                                 h_3 & g_2
                             \end{smallmatrix}
                           \right]\quad\;\;}&A_4\oplus B_3\oplus C_2\ar[r]^{\left[
                             \begin{smallmatrix}
                                -a_4 & 0 & 0\\
                                 -f_4 & -b_3 & 0\\
                                 h_4 & g_3 & c_2
                             \end{smallmatrix}
                           \right]\;\;}\ar[d]^{\left[
\begin{smallmatrix}
   u & v & w\\
      \end{smallmatrix}
 \right]}&A_4\oplus B_3\oplus C_2\ar@{.>}[dl]^{\hspace{5mm}\left[
  \begin{smallmatrix}
     l & m & n\\
      \end{smallmatrix}
\right]}\\
      &X&}$$
This shows that $\left[
                             \begin{matrix}
                                -a_4 & 0 & 0 \\
                                 -f_4 & -b_3 & 0 \\
                                 h_4 & g_3 & c_2 \\
                             \end{matrix}
                           \right]$ is a weak cokernel of $\left[
                             \begin{matrix}
                                -a_3 & 0 \\
                                 f_3 & -b_2 \\
                                 h_3 & g_2
                             \end{matrix}
                           \right]$.
Thus we prove that the sequence (\ref{RN4 last direct sums diagram induce}) is a right $n$-exact sequence.

We claim that the morphism $\left[
                    \begin{matrix}
                      a_2 \\
                      f_2 \\
                    \end{matrix}
                  \right]: A_2\rightarrow A_3\oplus B_2$
is $\X$-monic. In fact, for each morphism $s: A_2\rightarrow X$, where $X \in \X$, then $sa_1: A_1\rightarrow X$, since $f_1$ is $\X$-monic, there exist a morphism $t: B_1\rightarrow X$ such that $sa_1=tf_1$. It follows that $tf_1a_0=sa_1a_0=0$. So there exist a morphism $r: B_2\rightarrow X$ such that $rb_1=t$. Thus we have
$$(s-rf_2)a_1=sa_1-rf_2a_1=tf_1-rb_1f_1=tf_1-tf_1=0.$$
Then there exist a morphism $k: A_3\rightarrow X$ such that $ka_2=s-rf_2$.
It follows that $s=ka_2+rf_2=\left[
                    \begin{matrix}
                      k & r\\
                    \end{matrix}
                  \right]\left[
                    \begin{matrix}
                      a_2 \\
                      f_2 \\
                    \end{matrix}
                  \right]$.
This shows that
$\left[
                    \begin{matrix}
                      a_2 \\
                      f_2 \\
                    \end{matrix}
                  \right]: A_2\rightarrow A_3\oplus B_2$
is $\X$-monic. Hence the right $n$-exact sequence (\ref{RN4 last direct sums diagram induce}) induces a right $(n+2)$-angle
$$
A_2 \xrightarrow{\left[
                    \begin{smallmatrix}
                      \overline{a_2} \\[1mm]
                      \overline{f_2} \\
                    \end{smallmatrix}
                  \right]}A_3\oplus B_2 \xrightarrow{\left[
                             \begin{smallmatrix}
                                -\overline{a_3} & 0 \\[1mm]
                                 \overline{f_3} & -\overline{b_2} \\[1mm]
                                 \overline{h_3} & \overline{g_2}
                             \end{smallmatrix}
                           \right]}A_4 \oplus B_3\oplus C_2\xrightarrow{~\theta_1~}\cdots\hspace{50mm}$$
$$
\cdots \xrightarrow{~\theta_{n-3}~}A_{n+1}\oplus B_n\oplus C_{n-1} \xrightarrow{~\gamma~}B_{n+1}\oplus C_n \xrightarrow{\left[
                            \begin{smallmatrix}
                                 \overline{g_{n+1}} ~& \overline{c_n} \\
                             \end{smallmatrix}
                             \right]}C_{n+1}\xrightarrow{(-1)^n\overline{d_{n+1}}}\Sigma A_2,
$$
where $\theta_i=\left[\begin{matrix}
   -\overline{a_{i+3}} & 0 & 0\\[1mm]
  (-1)^i\overline{f_{i+3}} & -\overline{b_{i+2}} & 0\\[1mm]
   \overline{h_{i+3}} & \overline{g_{i+2}} & \overline{c_{i+1}}\\
 \end{matrix}\right],~\gamma=\left[\begin{matrix}
  (-1)^n\overline{f_{n+1}} & -\overline{b_{n}} & 0\\[1mm]
   \overline{h_{n+1}} & \overline{g_{n}} & \overline{c_{n-1}}\\
 \end{matrix}\right]$.\\[2mm]
For the following commutative diagram of right $n$-exact sequences
$$\xymatrix@C=2.5cm@R=1.5cm{A_1 \ar[r]^{f_1} \ar[d]^{a_1} & B_1 \ar[r]^{c_1} \ar[d]^{\left[
                             \begin{smallmatrix}
                                0 \\
                                b_1 \\
                             \end{smallmatrix}
                           \right]} & C_2 \ar[r]^{c_2} \ar[d]^{\left[
                             \begin{smallmatrix}
                                0 \\
                                0 \\
                                1 \\
                             \end{smallmatrix}
                           \right]} & \cdots\\
A_2 \ar[r]^{\left[
                    \begin{smallmatrix}
                      a_2 \\
                      f_2 \\
                    \end{smallmatrix}
                  \right]\hspace{4mm}} & A_3\oplus B_2 \ar[r]^{\left[
                             \begin{smallmatrix}
                               -a_3 & 0 \\
                               f_3 & -b_2 \\
                               h_3 & g_2 \\
                             \end{smallmatrix}
                           \right]\hspace{4mm}}& A_4\oplus B_3\oplus C_2 \ar[r]^{\hspace{8mm}\left[
                             \begin{smallmatrix}
                               -a_4 & 0 & 0\\
                               -f_4 & -b_3 &0\\
                               h_4 & g_3 & c_2\\
                             \end{smallmatrix}
                           \right]}& \cdots
}\hspace{45mm}$$
$$\hspace{45mm}\xymatrix@C=3.5cm@R=1.5cm{\cdots \ar[r]^{c_n-1}& C_n \ar[r]^{c_n} \ar[d]^{\left[
                             \begin{smallmatrix}
                               0 \\
                               1 \\
                             \end{smallmatrix}
                           \right]} & C_{n+1} \ar@{=}[d]\\
\cdots \ar[r]^{\left[
                             \begin{smallmatrix}
                              (-1)^{n}f_{n+1} & -b_n & 0 \\
                               h_{n+1} & g_n & c_{n-1} \\
                             \end{smallmatrix}
                           \right]\hspace{7mm}}& B_{n+1}\oplus C_n \ar[r]^{\hspace{2mm}\left[
\begin{smallmatrix}
                               g_{n+1} ~& c_{n} \\
                             \end{smallmatrix}
                           \right]}& C_{n+1}
}$$
by Lemma \ref{lemma 2.7 induced standard (n+2)-angle diagram}, we obtain the following commutative diagram of right $(n+2)$-angles
$$
\xymatrix@C=2.5cm@R=1.5cm{A_1 \ar[r]^{\overline{f_1}} \ar[d]^{\overline{a_1}} & B_1 \ar[r]^{\overline{c_1}} \ar[d]^{\left[
                             \begin{smallmatrix}
                                0 \\[1mm]
                                \overline{b_1} \\
                             \end{smallmatrix}
                           \right]} & C_2 \ar[r]^{\overline{c_2}} \ar[d]^{\left[
                             \begin{smallmatrix}
                                0 \\
                                0 \\
                                1 \\
                             \end{smallmatrix}
                           \right]} & \cdots\\
A_2 \ar[r]^{\left[
                    \begin{smallmatrix}
                      \overline{a_2} \\[1mm]
                      \overline{f_2} \\
                    \end{smallmatrix}
                  \right]\hspace{4mm}} & A_3\oplus B_2 \ar[r]^{\left[
                             \begin{smallmatrix}
                               -\overline{a_3} & 0 \\[1mm]
                               \overline{f_3} & -\overline{b_2} \\[1mm]
                               \overline{h_3} & \overline{g_2} \\
                             \end{smallmatrix}
                           \right]\hspace{4mm}}& A_4\oplus B_3\oplus C_2 \ar[r]^{\hspace{9mm}\left[
                             \begin{smallmatrix}
                               -\overline{a_4} & 0 & 0\\[1mm]
                               -\overline{f_4} & -\overline{b_3} &0\\[1mm]
                               \overline{h_4} & \overline{g_3} & \overline{c_2}\\
                             \end{smallmatrix}
                           \right]}& \cdots
}\hspace{45mm}$$
$$\xymatrix@C=3.5cm@R=1.5cm{\cdots \ar[r]^{\overline{c_{n-1}}}& C_n \ar[r]^{\overline{c_n}} \ar[d]^{\left[
                             \begin{smallmatrix}
                               0 \\
                               1 \\
                             \end{smallmatrix}
                           \right]} & C_{n+1} \ar[r]^{(-1)^n\overline{c_{n+1}}} \ar@{=}[d] & \Sigma A_1 \ar[d]^{\Sigma \overline{a_1}}\\
\cdots \ar[r]^{\left[
                             \begin{smallmatrix}
                              (-1)^{n}\overline{f_{n+1}} & -\overline{b_n} & 0 \\[1mm]
                               \overline{h_{n+1}} & \overline{g_n} & \overline{c_{n-1}} \\
                             \end{smallmatrix}
                           \right]\hspace{7mm}}& B_{n+1}\oplus C_n \ar[r]^{\left[
\begin{smallmatrix}
                               \overline{g_{n+1}} ~& \overline{c_{n}} \\
                             \end{smallmatrix}
                           \right]}& C_{n+1} \ar[r]^{(-1)^n\overline{d_{n+1}}} & \Sigma A_2
}$$
Thus we have $\overline{d_{n+1}}=\Sigma \overline{a_1}\circ \overline{c_{n+1}}$.
This shows that (RN4$^*$) holds.  \qed
\vspace{2mm}

In order to give an application, let's review some concepts from \cite{j}.

Let $\A$ be an $n$-abelian category. An object $I\in\A$ is called \emph{injective} if for any monomorphism $f:A\rightarrowtail B$, the sequence $$\A(B,I)\xrightarrow{\A(f,\hspace{1mm}I)}\A(A,I)\rightarrow 0$$ is exact. We denote by $\mathcal{I}$ the subcategory of injective objects. We say that $\A$ \emph{has enough injectives} if for any object  $A\in\A$, there exists an $n$-exact sequence
$$A\rightarrowtail I_1\rightarrow I_2\rightarrow\cdots\rightarrow I_n\twoheadrightarrow B$$
where $I_1,I_2,\cdots,I_n\in\mathcal{I}$.
\vspace{2mm}

As an application of Theorem \ref{main2}, we have the following conclusion.

\begin{corollary}
Let $\A$ be an $n$-abelian category with enough injectives and
$\mathcal{I}$ be the subcategory of injective objects. Then
the quotient category $\A/\mathcal{I}$ is a right $(n+2)$-angulated category.
\end{corollary}

\proof Since any morphism in $\A$ has $n$-cokernels and every left $\mathcal I$-approximation has a special $n$-cokernel with respect to $\mathcal I$,
this follows from Theorem \ref{main2}.  \qed
\vspace{2mm}

Now let's review the definition of a strongly covariantly finite
subcategory from \cite{LZ}.

\begin{definition}\label{dd1}\cite[Definition 3.1]{LZ}
Let $(\C,\Sigma,\Phi)$ be an $(n+2)$-angulated category. A subcategory $\X$ of $\C$ is called
\emph{strongly covariantly finite}, if for any object $B\in\C$, there exist an $(n+2)$-angle
$$B\xrightarrow{~f~}X_1\xrightarrow{}X_2\xrightarrow{}\cdots\xrightarrow{}X_{n}\xrightarrow{~}C\xrightarrow{~}\Sigma B$$
where $f$ is a left $\X$-approximation of $B$ and $X_1,X_2,\cdots,X_n\in\X$.
\end{definition}

Now we give some examples of strongly covariance finite
subcategories.

\begin{example}
Let $(\C,\Sigma,\Phi)$ be an $(n+2)$-angulated category and $\X$ be a cluster tilting subcategory (in the sense of Zhou and Zhu \cite[Definition 1.1]{ZZ}) of $\C$.
By Definition \ref{dd1}, we know that $\X$ is strongly covariantly finite of $\C$.
\end{example}

\begin{example}
This example comes from \cite{L1}.
Let $$\T=D^b(kQ)/\tau^{-1}[1]$$ be the cluster category of type $A_3$,
where $Q$ is the quiver $1\xrightarrow{~\alpha~}2\xrightarrow{~\beta~}3$,
$D^b(kQ)$ is the bounded derived category of finite generated modules over $kQ$, $\tau$ is the Auslander-Reiten translation and $[1]$ is the shift functor
of $D^b(kQ)$. Then $\T$ is a triangulated category. Its shift functor is also denoted by $[1]$.
We describe the
Auslander-Reiten quiver of $\T$ in the following:
$$\xymatrix@C=0.6cm@R0.3cm{
&&P_1\ar[dr]
&&S_3[1]\ar[dr]
&&\\
&P_2 \ar@{.}[rr] \ar[dr] \ar[ur]
&&I_2 \ar@{.}[rr] \ar[dr] \ar[ur]
&&P_2[1]\ar[dr]\\
S_3\ar[ur]\ar@{.}[rr]&&S_2\ar[ur]\ar@{.}[rr]
&&S_1\ar[ur]\ar@{.}[rr]&&P_1[1]
}
$$
It is straightforward to verify that $\C:=\add(S_3\oplus P_1\oplus S_1)$ is a $2$-cluster tilting subcategory of $\T$ and satisfies $\C[2]=\C$.   By \cite[Theorem 1]{GKO}, we know that $\C$ is a $4$-angulated category with an automorphism functor $[2]$.
Let $\X=\add(S_3\oplus S_1)$.
Then the $4$-angle
$$P_1\xrightarrow{~~}S_1\xrightarrow{~~}S_3\xrightarrow{~~}P_1\xrightarrow{~~}P_1[2]$$
shows that $\X$ is a strongly covariantly finite subcategory of $\C$.
\end{example}

Now we provide another example of a right $(n+2)$-angulated category.

\begin{definition}{\rm \cite[Proposition 3.3]{L1}}\label{definition 4.12 functor}
Let $(\C,\Sigma,\Phi)$ be an $(n+2)$-angulated category and $\X$ be a strongly covariantly finite subcategory of $\C$. Then there exist an additive endofunctor $\mathbb{G}\colon\C / \X \to \C / \X$, defined as follows:

For any object $A\in \C$, there exists an $(n+2)$-angle
$$A\xrightarrow{~a_0~}X_1\xrightarrow{~a_1~}X_2\xrightarrow{~a_2~}\cdots \xrightarrow{~a_{n-1}~} X_n\xrightarrow{~a_n~}B\xrightarrow{~a_{n+1}~}\Sigma A$$
where $X_1, \cdots, X_n \in \X$, $a_0$ is a left $\X$-approximation of $A$. For any morphism $f\colon A\rightarrow A'$, since $a_0$ is a left $\X$-approximation of $A$, we have the following commutative diagram
$$\xymatrix{
A\ar[r]^{a_0}\ar[d]^{f} & X_1\ar[r]^{a_1}\ar[d]^{f_1} & X_2\ar[r]^{a_2}\ar[d]^{f_2} & \cdots \ar[r]^{a_{n-1}} & X_{n}\ar[r]^{a_{n}}\ar[d]^{f_{n}}  & B \ar[d]^{g} \ar[r]^{a_{n+1}} & \Sigma A \ar[d]^{\Sigma f}\\
A'\ar[r]^{a_0'} & X'_1\ar[r]^{a_1'} & X'_2\ar[r]^{a_2'} & \cdots \ar[r]^{a_{n-1}'} & X'_{n}\ar[r]^{a_{n}'} &  B'\ar[r]^{a_{n+1}'} & \Sigma A'.
}$$
Define an endofunctor $\mathbb{G}\colon\C/\X\rightarrow\C/\X$ such that $\mathbb{G}A=B$ and $\mathbb{G} \overline{f}=\overline{g}$.
\end{definition}

\begin{definition}\cite[Definition 3.4]{L1}\label{definition 4.13 standard (n+2)-angle}
Let $(\C,\Sigma,\Phi)$ be an $(n+2)$-angulated category and $\X$ be a strongly covariantly finite subcategory of $\C$. Assume that
$$\xymatrix{A_0 \ar[r]^{a_0} & A_1 \ar[r]^{a_1} & A_2 \ar[r]^{a_2} & \cdots \ar[r]^{a_{n-1}} & A_n \ar[r]^{a_n\hspace{2mm}} & A_{n+1} \ar[r]^{a_{n+1}} & \Sigma A_0}$$
is an $(n+2)$-angle in $\C$, where $a_0$ is $\X$-monic. Then there exist a commutative diagram
$$\xymatrix{
A_0 \ar[r]^{a_0}\ar@{=}[d] & A_1 \ar[r]^{a_1}\ar[d]^{f_1} & A_2 \ar[r]^{a_2}\ar[d]^{f_2} & \cdots \ar[r]^{a_{n-1}}& A_{n} \ar[r]^{a_{n}\hspace{2mm}}\ar[d]^{f_{n}} & A_{n+1} \ar[d]^{f_{n+1}} \ar[r]^{a_{n+1}} & \Sigma A_0 \ar@{=}[d]\\
A_0 \ar[r]^{b_0} & X_1 \ar[r]^{b_1} & X_2 \ar[r]^{b_2} & \cdots \ar[r]^{b_{n-1}}& X_{n} \ar[r]^{b_{n}\hspace{2mm}} & \mathbb{G}A_0 \ar[r]^{b_{n+1}} & \Sigma A_0
}$$
of $(n+2)$-angles. Then we have a complex
$$A_0 \xrightarrow{~\overline{a_0}~} A_1\xrightarrow{~\overline{a_1}~} A_2 \xrightarrow{~\overline{a_2}~} \cdots \xrightarrow{~\overline{a_n}~} A_{n+1} \xrightarrow{(-1)^n\overline{f_{n+1}}} \mathbb{G}A_0.
$$
We define \emph{right $(n+2)$-angles} in $\C/\X$ as the complexes which are isomorphic to
complexes obtained in this way. We denote by $\Theta$ the class of right $(n+2)$-angles.
\end{definition}

The following result was proved in \cite[Lemma 3.5]{L1} for general case.
So their proof can be applied for our case without any change.

\begin{lemma}\emph{\cite[Lemma 3.5]{L1}}\label{lemma 4.14 induced standard (n+2)-angle diagram 2}
Let
$$\xymatrix{A_0 \ar[r]^{a_0} \ar[d]^{f_0} & A_1 \ar[r]^{a_1} \ar[d]^{f_1} & A_2 \ar[r]^{a_2} \ar[d]^{f_2} &\cdots \ar[r]^{a_{n-1}} & A_n \ar[r]^{a_n\hspace{2mm}} \ar[d]^{f_n} & A_{n+1} \ar[d]^{f_{n+1}} \ar[r]^{a_{n+1}} & \Sigma A_0 \ar[d]^{\Sigma f_0}\\
B_0 \ar[r]^{b_0} & B_1 \ar[r]^{b_1} & B_2 \ar[r]^{b_2} &\cdots \ar[r]^{b_{n-1}} & B_n \ar[r]^{b_n\hspace{2mm}} & B_{n+1} \ar[r]^{b_{n+1}} & \Sigma B_0
}$$
be a commutative diagram of $(n+2)$-angles in $\C$, where $a_0$ and $b_0$ are $\X$-monic. Then we have a commutative diagram
$$\xymatrix@C=1.5cm{A_0 \ar[r]^{\overline{a_0}} \ar[d]^{\overline{f_0}} & A_1 \ar[r]^{\overline{a_1}} \ar[d]^{\overline{f_1}} & A_2 \ar[r]^{\overline{a_2}} \ar[d]^{\overline{f_2}} &\cdots \ar[r]^{\overline{a_{n-1}}} & A_n \ar[r]^{\overline{a_n}\hspace{2mm}} \ar[d]^{\overline{f_n}} & A_{n+1} \ar[r]^{(-1)^n\overline{\alpha_{n+1}}} \ar[d]^{\overline{f_{n+1}}} & \mathbb{G}A_0 \ar[d]^{\mathbb{G}\overline{f_0}}\\
B_0 \ar[r]^{\overline{b_0}} & B_1 \ar[r]^{\overline{b_1}} & B_2 \ar[r]^{\overline{b_2}} &\cdots \ar[r]^{\overline{b_{n-1}}} & B_n \ar[r]^{\overline{b_n}\hspace{2mm}} & B_{n+1} \ar[r]^{(-1)^n\overline{\beta_{n+1}}} & \mathbb{G}B_0
}$$
of right $(n+2)$-angles in $\C/\X$.
\end{lemma}

\begin{theorem}
Let $(\C,\Sigma,\Phi)$ be an $(n+2)$-angulated category and $\X$ be a
strongly covariantly finite subcategory of $\C$.
Then the quotient category $\C/\X$ is a right $(n+2)$-angulated category
in the sense of Definition \ref{right (n+2)-angulated category}, where the endofunctor $\mathbb{G}$ is defined in Definition \emph{\ref{definition 4.12 functor}} and right $(n+2)$-angles is defined in Definition \emph{\ref{definition 4.13 standard (n+2)-angle}}.
\end{theorem}

\proof We only need to prove (RN1)(b$^*$) and (RN4$^*$), as the proofs for the other axioms are the same as the ones provided in \cite[Theorem 3.7]{L1} and can be omitted.

 The commutative diagram
$$\xymatrix{0 \ar[r] \ar@{=}[d] & A \ar[r]^{1} \ar[d] & A \ar[r] \ar[d] & 0 \ar[r] \ar[d] & \cdots \ar[r] & 0 \ar[r] \ar[d] &0 \ar@{=}[d]\\
0 \ar[r] & 0 \ar[r] & 0 \ar[r] &0 \ar[r] & \cdots \ar[r] & 0 \ar[r] & 0
}$$
shows that
$0\rightarrow A\xrightarrow{1} A\rightarrow 0\rightarrow \cdots\rightarrow 0$
belongs to $\Theta$. Thus (RN1)(b$^*$) holds.
Now we prove that (RN4$^*$) holds. Given the solid part of the diagram
$$
\xymatrix@C=1.5cm{A_0 \ar[r]^{\overline{a_0}} \ar@{=}[d] & A_1 \ar[r]^{\overline{a_1}} \ar[d]^{\overline{f_1}} & A_2 \ar[r]^{\overline{a_2}} &\cdots  \ar[r]^{\overline{a_n}} & A_{n+1} \ar[r]^{(-1)^n\overline{\alpha_{n+1}}} & \mathbb{G}A_0 \ar@{=}[d]\\
A_0 \ar[r]^{\overline{f_1a_0}} \ar[d]^{\overline{a_0}} & B_1 \ar[r]^{\overline{b_1}} \ar@{=}[d] & B_2 \ar[r]^{\overline{b_2}} &\cdots \ar[r]^{\overline{b_n}} & B_{n+1} \ar[r]^{(-1)^n\overline{\beta_{n+1}}} & \mathbb{G}A_0 \ar[d]^{\mathbb{G}\overline{a_0}}\\
A_1 \ar[r]^{\overline{f_1}} & B_1 \ar[r]^{\overline{c_1}} & C_2 \ar[r]^{\overline{c_2}} &\cdots \ar[r]^{\overline{c_n}} & C_{n+1} \ar[r]^{(-1)^n\overline{\theta_{n+1}}} & \mathbb{G}A_1
}
$$
with commuting squares and with rows in $\Theta$ which are induced by the three $(n+2)$-angles in $\C$ as follows
$$A_0 \xrightarrow{~a_0~} A_1\xrightarrow{~a_1~} A_2 \xrightarrow{~a_2~} \cdots \xrightarrow{~a_n~} A_{n+1} \xrightarrow{a_{n+1}} \Sigma A_0,
$$
$$A_0 \xrightarrow{~f_1a_0~} B_1\xrightarrow{~b_1~} B_2 \xrightarrow{~b_2~} \cdots \xrightarrow{~b_n~} B_{n+1} \xrightarrow{b_{n+1}} \Sigma A_0,
$$
$$A_1 \xrightarrow{~f_1~} B_1\xrightarrow{~c_1~} C_2 \xrightarrow{~c_2~} \cdots \xrightarrow{~c_n~} C_{n+1} \xrightarrow{c_{n+1}} \Sigma A_1.
$$
where $a_0$ and $f_1$ are $\X$-monic, so is $f_1a_0$. Apply (N4$^*$) to the above three $(n+2)$-angles, there exist morphisms
$
f_i: A_i \rightarrow B_i~~(i=2, 3, \cdots, n+1),
g_i: B_i \rightarrow C_i~~(i=2, 3, \cdots, n+1),
h_i: A_i \rightarrow C_{i-1}~~(i=3, 4, \cdots, n+1)$
such that each square of the following diagram are commutative
\begin{equation}\label{th 4.14 n4^* diagram}
\begin{split}
\xymatrix{A_0 \ar[r]^{a_0} \ar@{=}[d] & A_1 \ar[r]^{a_1} \ar[d]^{f_1} & A_2 \ar[r]^{a_2} \ar@{-->}[d]^{f_2} &\cdots \ar[r]^{a_{n-1}} & A_n \ar[r]^{a_n\hspace{2mm}} \ar@{-->}[d]^{f_n} & A_{n+1} \ar[r]^{a_{n+1}} \ar@{-->}[d]^{f_{n+1}} & \Sigma A_0 \ar@{=}[d]\\
A_0 \ar[r]^{f_1a_0} \ar[d]^{a_0} & B_1 \ar[r]^{b_1} \ar@{=}[d] & B_2 \ar[r]^{b_2} \ar@{-->}[d]^{g_2} &\cdots \ar[r]^{b_{n-1}} & B_n \ar[r]^{b_n\hspace{2mm}} \ar@{-->}[d]^{g_n} & B_{n+1} \ar[r]^{b_{n+1}} \ar@{-->}[d]^{g_{n+1}} & \Sigma A_0 \ar[d]^{\Sigma a_0}\\
A_1 \ar[r]^{f_1} & B_1 \ar[r]^{c_1} & C_2 \ar[r]^{c_2} &\cdots \ar[r]^{c_{n-1}} & C_n \ar[r]^{c_n\hspace{2mm}} & C_{n+1} \ar[r]^{c_{n+1}} & \Sigma A_1
}
\end{split}
\end{equation}
and the following $(n+2)$-$\Sigma$-sequence
$$A_2 \xrightarrow{\left[
                    \begin{smallmatrix}
                      a_2 \\
                      f_2 \\
                    \end{smallmatrix}
                  \right]}A_3\oplus B_2 \xrightarrow{\left[
                             \begin{smallmatrix}
                                -a_3 & 0 \\
                                 f_3 & -b_2 \\
                                 h_3 & g_2
                             \end{smallmatrix}
                           \right]}A_4 \oplus B_3\oplus C_2\xrightarrow{\left[
                             \begin{smallmatrix}
                                -a_4 & 0 & 0 \\
                                 -f_4 & -b_3 & 0 \\
                                 h_4 & g_3 & c_2 \\
                             \end{smallmatrix}
                           \right]}\cdots\hspace{30mm}$$
\begin{equation}\label{th 4.14 apply n4^* sequence}
\cdots \xrightarrow{\left[
 \begin{smallmatrix}
     (-1)^nf_{n+1} & -b_n & 0 \\
   h_{n+1} & g_n & c_{n-1} \\
  \end{smallmatrix}
  \right]}B_{n+1}\oplus C_n \xrightarrow{\left[
                            \begin{smallmatrix}
                                 g_{n+1} & c_n \\
                             \end{smallmatrix}
                             \right]}C_{n+1} \xrightarrow{\Sigma a_1 \circ c_{n+1}}\Sigma A_2
\end{equation}
belongs to $\Phi$.
From the diagram (\ref{th 4.14 n4^* diagram}) and Lemma \ref{lemma 4.14 induced standard (n+2)-angle diagram 2}, we get the following commutative diagram of right $(n+2)$-angles
$$
\xymatrix@C=1.5cm{A_0 \ar[r]^{\overline{a_0}} \ar@{=}[d] & A_1 \ar[r]^{\overline{a_1}} \ar[d]^{\overline{f_1}} & A_2 \ar[r]^{\overline{a_2}} \ar[d]^{\overline{f_2}} &\cdots \ar[r]^{\overline{a_{n-1}}} & A_n \ar[r]^{\overline{a_n}\hspace{2mm}} \ar[d]^{\overline{f_n}} & A_{n+1} \ar[r]^{(-1)^n\overline{\alpha_{n+1}}} \ar[d]^{\overline{f_{n+1}}} & \mathbb{G}A_0 \ar@{=}[d]\\
A_0 \ar[r]^{\overline{{f_1a_0}}} \ar[d]^{\overline{a_0}} & B_1 \ar[r]^{\overline{b_1}} \ar@{=}[d] & B_2 \ar[r]^{\overline{b_2}} \ar[d]^{\overline{g_2}} &\cdots \ar[r]^{\overline{b_{n-1}}} & B_n \ar[r]^{\overline{b_n}\hspace{2mm}} \ar[d]^{\overline{g_n}} & B_{n+1} \ar[r]^{(-1)^n\overline{\beta_{n+1}}} \ar[d]^{\overline{g_{n+1}}} & \mathbb{G}A_0 \ar[d]^{\mathbb{G}\overline{a_0}}\\
A_1 \ar[r]^{\overline{f_1}} & B_1 \ar[r]^{\overline{c_1}} & C_2 \ar[r]^{\overline{c_2}} &\cdots \ar[r]^{\overline{c_{n-1}}} & C_n \ar[r]^{\overline{c_n}\hspace{2mm}} & C_{n+1} \ar[r]^{(-1)^n\overline{\theta_{n+1}}} & \mathbb{G}A_1
}
$$
We claim that the morphism $\left[
                    \begin{matrix}
                      a_2 \\
                      f_2 \\
                    \end{matrix}
                  \right]: A_2\rightarrow A_3\oplus B_2$
is $\X$-monic. In fact, for each morphism $s: A_2\rightarrow X$, where $X \in \X$, since $f_1$ is $\X$-monic, there exist a morphism $t: B_1\rightarrow X$ such that $sa_1=tf_1$. It follows that $tf_1a_0=sa_1a_0=0$. So there exist a morphism $r: B_2\rightarrow X$ such that $rb_1=t$. Thus we have
$$(s-rf_2)a_1=sa_1-rf_2a_1=tf_1-rb_1f_1=tf_1-tf_1=0.$$
Then there exist a morphism $k: A_3\rightarrow X$ such that $ka_2=s-rf_2$.
It follows that $s=ka_2+rf_2=\left[
                    \begin{matrix}
                      k & r\\
                    \end{matrix}
                  \right]\left[
                    \begin{matrix}
                      a_2 \\
                      f_2 \\
                    \end{matrix}
                  \right]$.
This shows that
$\left[
                    \begin{matrix}
                      a_2 \\
                      f_2 \\
                    \end{matrix}
                  \right]: A_2\rightarrow A_3\oplus B_2$
is $\X$-monic. Hence the sequence (\ref{th 4.14 apply n4^* sequence}) induces a right $(n+2)$-angle
$$A_2 \xrightarrow{\left[
                    \begin{smallmatrix}
                      \overline{a_2} \\[1mm]
                      \overline{f_2} \\
                    \end{smallmatrix}
                  \right]}A_3\oplus B_2 \xrightarrow{\left[
                             \begin{smallmatrix}
                                -\overline{a_3} & 0 \\[1mm]
                                 \overline{f_3} & -\overline{b_2} \\[1mm]
                                 \overline{h_3} & \overline{g_2}
                             \end{smallmatrix}
                           \right]}A_4 \oplus B_3\oplus C_2\xrightarrow{\left[
                             \begin{smallmatrix}
                                -\overline{a_4} & 0 & 0 \\[1mm]
                               -\overline{f_4}& -\overline{b_3} & 0 \\[1mm]
                             \overline{h_4} & \overline{g_3} & \overline{c_2} \\
                             \end{smallmatrix}
                           \right]}\cdots\hspace{30mm}$$
$$
\cdots \xrightarrow{\left[
 \begin{smallmatrix}
     (-1)^n\overline{f_{n+1}} & -\overline{b_n} & 0 \\[1mm]
   \overline{h_{n+1}} & \overline{g_n} & \overline{c_{n-1}} \\
  \end{smallmatrix}
  \right]}B_{n+1}\oplus C_n \xrightarrow{\left[
                            \begin{smallmatrix}
                                 \overline{g_{n+1}} & \overline{c_n} \\
                             \end{smallmatrix}
                             \right]}C_{n+1} \xrightarrow{(-1)^n\overline{\gamma_{n+1}} }\mathbb{G}A_2.
$$
For the commutative diagram of $(n+2)$-angles
$$\xymatrix@C=2.5cm@R=1.5cm{A_1 \ar[r]^{f_1} \ar[d]^{a_1} & B_1 \ar[r]^{c_1} \ar[d]^{\left[
                             \begin{smallmatrix}
                                0 \\
                                b_1 \\
                             \end{smallmatrix}
                           \right]} & C_2 \ar[r]^{c_2} \ar[d]^{\left[
                             \begin{smallmatrix}
                                0 \\
                                0 \\
                                1 \\
                             \end{smallmatrix}
                           \right]} & \cdots\\
A_2 \ar[r]^{\left[
                    \begin{smallmatrix}
                      a_2 \\
                      f_2 \\
                    \end{smallmatrix}
                  \right]\hspace{4mm}} & A_3\oplus B_2 \ar[r]^{\left[
                             \begin{smallmatrix}
                               -a_3 & 0 \\
                               f_3 & -b_2 \\
                               h_3 & g_2 \\
                             \end{smallmatrix}
                           \right]\hspace{4mm}}& A_4\oplus B_3\oplus C_2 \ar[r]^{\hspace{8mm}\left[
                             \begin{smallmatrix}
                               -a_4 & 0 & 0\\
                               -f_4 & -b_3 &0\\
                               h_4 & g_3 & c_2\\
                             \end{smallmatrix}
                           \right]}& \cdots
}\hspace{45mm}$$
$$\xymatrix@C=3.5cm@R=1.5cm{\cdots \ar[r]^{c_n-1}& C_n \ar[r]^{c_n} \ar[d]^{\left[
                             \begin{smallmatrix}
                               0 \\
                               1 \\
                             \end{smallmatrix}
                           \right]} & C_{n+1} \ar@{=}[d] \ar[r]^{c_{n+1}} & \Sigma A_1 \ar[d]^{\Sigma a_1}\\
\cdots \ar[r]^{\left[
                             \begin{smallmatrix}
                              (-1)^{n}f_{n+1} & -b_n & 0 \\
                               h_{n+1} & g_n & c_{n-1} \\
                             \end{smallmatrix}
                           \right]\hspace{7mm}}& B_{n+1}\oplus C_n \ar[r]^{\hspace{2mm}\left[
\begin{smallmatrix}
                               g_{n+1} ~& c_{n} \\
                             \end{smallmatrix}
                           \right]}& C_{n+1} \ar[r]^{\Sigma a_1\circ c_{n+1}} &\Sigma A_2
}$$
by Lemma \ref{lemma 4.14 induced standard (n+2)-angle diagram 2}, we obtain the following commutative diagram of right $(n+2)$-angles
$$
\xymatrix@C=2.5cm@R=1.5cm{A_1 \ar[r]^{\overline{f_1}} \ar[d]^{\overline{a_1}} & B_1 \ar[r]^{\overline{c_1}} \ar[d]^{\left[
                             \begin{smallmatrix}
                                0 \\[1mm]
                                \overline{b_1} \\
                             \end{smallmatrix}
                           \right]} & C_2 \ar[r]^{\overline{c_2}} \ar[d]^{\left[
                             \begin{smallmatrix}
                                0 \\
                                0 \\
                                1 \\
                             \end{smallmatrix}
                           \right]} & \cdots\\
A_2 \ar[r]^{\left[
                    \begin{smallmatrix}
                      \overline{a_2} \\[1mm]
                      \overline{f_2} \\
                    \end{smallmatrix}
                  \right]\hspace{4mm}} & A_3\oplus B_2 \ar[r]^{\left[
                             \begin{smallmatrix}
                               -\overline{a_3} & 0 \\[1mm]
                               \overline{f_3} & -\overline{b_2} \\[1mm]
                               \overline{h_3} & \overline{g_2} \\
                             \end{smallmatrix}
                           \right]\hspace{4mm}}& A_4\oplus B_3\oplus C_2 \ar[r]^{\hspace{9mm}\left[
                             \begin{smallmatrix}
                               -\overline{a_4} & 0 & 0\\[1mm]
                               -\overline{f_4} & -\overline{b_3} &0\\[1mm]
                               \overline{h_4} & \overline{g_3} & \overline{c_2}\\
                             \end{smallmatrix}
                           \right]}& \cdots
}\hspace{45mm}$$
$$\xymatrix@C=3.5cm@R=1.5cm{\cdots \ar[r]^{\overline{c_{n-1}}}& C_n \ar[r]^{\overline{c_n}} \ar[d]^{\left[
                             \begin{smallmatrix}
                               0 \\
                               1 \\
                             \end{smallmatrix}
                           \right]} & C_{n+1} \ar[r]^{(-1)^n\overline{\theta_{n+1}}} \ar@{=}[d] & \mathbb{G}A_1 \ar[d]^{\mathbb{G}\overline{a_1}}\\
\cdots \ar[r]^{\left[
                             \begin{smallmatrix}
                              (-1)^{n}\overline{f_{n+1}} & -\overline{b_n} & 0 \\[1mm]
                               \overline{h_{n+1}} & \overline{g_n} & \overline{c_{n-1}} \\
                             \end{smallmatrix}
                           \right]\hspace{7mm}}& B_{n+1}\oplus C_n \ar[r]^{\left[
\begin{smallmatrix}
                               \overline{g_{n+1}} ~& \overline{c_{n}} \\
                             \end{smallmatrix}
                           \right]}& C_{n+1} \ar[r]^{(-1)^n\overline{\gamma_{n+1}}} & \mathbb{G}A_2
}$$
Thus we have $\overline{\gamma_{n+1}}=\mathbb{G}\overline{a_1}\circ \overline{\theta_{n+1}}$.
This shows that (RN4$^*$) holds.\qed

\section{Axiom (RN3) is redundant}
\setcounter{equation}{0}

In this section, we will prove that (RN3) is redundant, as it can be deduced from other axioms. In fact, it can be deduced by (RN1)(c) and (RN4$^\ast$).

\begin{theorem}\label{main1}
Let $(\C,\Sigma,\Theta)$ be a right $(n+2)$-angulated category.
Then the axiom \emph{(RN3)} is a consequence of the axioms \emph{(RN1)(c)} and \emph{(RN4$^*$)}.
\end{theorem}

\proof  This is an adaptation of the proof of \cite[Theorem 3.1]{ahbt}.
Given the solid part of the commutative diagram
$$\xymatrix{
A_0 \ar[r]^{a_0}\ar[d]^{f_0} & A_1 \ar[r]^{a_1}\ar[d]^{f_1} & A_2 \ar[r]^{a_2} & \cdots \ar[r]^{a_{n-1}}& A_{n} \ar[r]^{a_{n}\hspace{2mm}} & A_{n+1} \ar[r]^{a_{n+1}} & \Sigma A_0 \ar[d]^{\Sigma f_0}\\
B_0 \ar[r]^{b_0} & B_1 \ar[r]^{b_1} & B_2 \ar[r]^{b_2} & \cdots \ar[r]^{b_{n-1}}& B_{n} \ar[r]^{b_{n}\hspace{2mm}} & B_{n+1} \ar[r]^{b_{n+1}} & \Sigma B_0
}$$
with rows in $\Theta$. From the above diagram we have $f_1a_0=b_0f_0$, we write it as $c_0$, i.e. $$c_0=f_1a_0=b_0f_0.$$
By (RN1)(c), the three morphisms $c_0, f_0, f_1$ can be respectively embedded into three right $(n+2)$-angles
$$A_0 \xrightarrow{~c_0~} B_1 \xrightarrow{~c_1~} C_2 \xrightarrow{~c_2~} \cdots \xrightarrow{~c_{n-1}~}C_n \xrightarrow{~c_{n}~} C_{n+1} \xrightarrow{~c_{n+1}~} \Sigma A_0,$$
$$A_0 \xrightarrow{~f_0~} B_0 \xrightarrow{~d_1~} D_2 \xrightarrow{~d_2~} \cdots \xrightarrow{~d_{n-1}~}D_n \xrightarrow{~d_{n}~} D_{n+1} \xrightarrow{~d_{n+1}~} \Sigma A_0,$$
$$A_1 \xrightarrow{~f_1~} B_1 \xrightarrow{~e_1~} E_2 \xrightarrow{~e_2~} \cdots \xrightarrow{~e_{n-1}~}E_n \xrightarrow{~e_{n}~} E_{n+1} \xrightarrow{~e_{n+1}~} \Sigma A_1.$$
Consider the following two commutative diagram
\begin{equation}\label{RN4 1 diagram}
\begin{split}
\xymatrix{A_0 \ar[r]^{a_0} \ar@{=}[d] & A_1 \ar[r]^{a_1} \ar[d]^{f_1} & A_2 \ar[r]^{a_2} \ar@{-->}[d]^{\alpha_2} &\cdots \ar[r]^{a_{n-1}} & A_n \ar[r]^{a_n\hspace{2mm}} \ar@{-->}[d]^{\alpha_n} & A_{n+1} \ar[r]^{a_{n+1}} \ar@{-->}[d]^{\alpha_{n+1}} & \Sigma A_0 \ar@{=}[d]\\
A_0 \ar[r]^{c_0} \ar[d]^{a_0} & B_1 \ar[r]^{c_1} \ar@{=}[d] & C_2 \ar[r]^{c_2} \ar@{-->}[d]^{\beta_2} &\cdots \ar[r]^{c_{n-1}} & C_n \ar[r]^{c_n\hspace{2mm}} \ar@{-->}[d]^{\beta_n} & C_{n+1} \ar[r]^{c_{n+1}} \ar@{-->}[d]^{\beta_{n+1}} & \Sigma A_0 \ar[d]^{\Sigma a_0}\\
A_1 \ar[r]^{f_1} & B_1 \ar[r]^{e_1} & E_2 \ar[r]^{e_2} &\cdots \ar[r]^{e_{n-1}} & E_n \ar[r]^{e_n\hspace{2mm}} & E_{n+1} \ar[r]^{e_{n+1}} & \Sigma A_1
}
\end{split}
\end{equation}
\begin{equation}\label{RN4 2 diagram}
\begin{split}
\xymatrix{A_0 \ar[r]^{f_0} \ar@{=}[d] & B_0 \ar[r]^{d_1} \ar[d]^{b_0} & D_2 \ar[r]^{d_2} \ar@{-->}[d]^{\theta_2} &\cdots \ar[r]^{d_{n-1}} & D_n \ar[r]^{d_n\hspace{2mm}} \ar@{-->}[d]^{\theta_n} & D_{n+1} \ar[r]^{d_{n+1}} \ar@{-->}[d]^{\theta_{n+1}} & \Sigma A_0 \ar@{=}[d]\\
A_0 \ar[r]^{c_0} \ar[d]^{f_0} & B_1 \ar[r]^{c_1} \ar@{=}[d] & C_2 \ar[r]^{c_2} \ar@{-->}[d]^{\mu_2} &\cdots \ar[r]^{c_{n-1}} & C_n \ar[r]^{c_n\hspace{2mm}} \ar@{-->}[d]^{\mu_n} & C_{n+1} \ar[r]^{c_{n+1}} \ar@{-->}[d]^{\mu_{n+1}} & \Sigma A_0 \ar[d]^{\Sigma f_0}\\
B_0 \ar[r]^{b_0} & B_1 \ar[r]^{b_1} & B_2 \ar[r]^{b_2} &\cdots \ar[r]^{b_{n-1}} & B_n \ar[r]^{b_n\hspace{2mm}} & B_{n+1} \ar[r]^{b_{n+1}} & \Sigma B_0
}
\end{split}
\end{equation}
with rows in $\Theta$, for the two solid part of (\ref{RN4 1 diagram}) and (\ref{RN4 2 diagram}), applies (RN4$^*$), there exist dotted morphisms such that the two diagrams above commute.

Since the diagrams (\ref{RN4 1 diagram}) and (\ref{RN4 2 diagram}) have the same middle row. We can combine the top half of (\ref{RN4 1 diagram}) with the bottom half of (\ref{RN4 2 diagram}) to obtain the following commutative diagram
\begin{equation}\label{combine diagram}
\begin{split}
\xymatrix{A_0 \ar[r]^{a_0} \ar@{=}[d] & A_1 \ar[r]^{a_1} \ar[d]^{f_1} & A_2 \ar[r]^{a_2} \ar@{-->}[d]^{\alpha_2} &\cdots \ar[r]^{a_{n-1}} & A_n \ar[r]^{a_n\hspace{2mm}} \ar@{-->}[d]^{\alpha_n} & A_{n+1} \ar[r]^{a_{n+1}} \ar@{-->}[d]^{\alpha_{n+1}} & \Sigma A_0 \ar@{=}[d]\\
A_0 \ar[r]^{c_0} \ar[d]^{f_0} & B_1 \ar[r]^{c_1} \ar@{=}[d] & C_2 \ar[r]^{c_2} \ar@{-->}[d]^{\mu_2} &\cdots \ar[r]^{c_{n-1}} & C_n \ar[r]^{c_n\hspace{2mm}} \ar@{-->}[d]^{\mu_n} & C_{n+1} \ar[r]^{c_{n+1}} \ar@{-->}[d]^{\mu_{n+1}} & \Sigma A_0 \ar[d]^{\Sigma f_0}\\
B_0 \ar[r]^{b_0} & B_1 \ar[r]^{b_1} & B_2 \ar[r]^{b_2} &\cdots \ar[r]^{b_{n-1}} & B_n \ar[r]^{b_n\hspace{2mm}} & B_{n+1} \ar[r]^{b_{n+1}} & \Sigma B_0.
}
\end{split}
\end{equation}
By omitting the middle line of (\ref{combine diagram}), we obtain the following commutative diagram
$$\xymatrix{
A_0 \ar[r]^{a_0}\ar[d]^{f_0} & A_1 \ar[r]^{a_1}\ar[d]^{f_1} & A_2 \ar[r]^{a_2} \ar[d]^{\mu_2\alpha_2} & \cdots \ar[r]^{a_{n-1}}& A_{n} \ar[r]^{a_{n}\hspace{2mm}} \ar[d]^{\mu_n\alpha_n} & A_{n+1} \ar[r]^{a_{n+1}} \ar[d]^{\mu_{n+1}\alpha_{n+1}} & \Sigma A_0 \ar[d]^{\Sigma f_0}\\
B_0 \ar[r]^{b_0} & B_1 \ar[r]^{b_1} & B_2 \ar[r]^{b_2} & \cdots \ar[r]^{b_{n-1}}& B_{n} \ar[r]^{b_{n}\hspace{2mm}} & B_{n+1} \ar[r]^{b_{n+1}} & \Sigma B_0
}$$
We define $f_k= \mu_k\alpha_k$ for $2\leq k\leq n+1$, then $(f_0, f_1, f_2, \cdots, f_{n+1})$ is a morphism of right $(n+2)$-angles. Thus (RN3) holds.  \qed
\vspace{2mm}

This theorem immediately yields the following conclusion.

\begin{corollary}{\rm \cite[Theorem 3.1]{ahbt}}
Let $(\C,\Sigma,\Phi)$ be an $(n+2)$-angulated category.
Then the axiom \emph{(N3)} is a consequence of the axioms \emph{(N1)(c)} and \emph{(N4)}.
\end{corollary}

\proof Since any $(n+2)$-angulated category can be viewed as a right $(n+2)$-angulated category, this follows from Theorem \ref{main1}.  \qed

\section{Some equivalent characterizations of axiom (RN4$^*$)}
\setcounter{equation}{0}

In this section,  we will give some equivalent characterizations of axiom (RN4$^*$).

\begin{theorem}
If $\Theta$ is a collection of $(n+2)$-$\Sigma$-sequences satisfying the axioms \emph{(RN1)}, {\rm (RN2)} and {\rm (RN3)}, then the following statements are equivalent:

$(1)$ $\Theta$ satisfies \emph{(RN4$^*$);}

$(2)$ $\Theta$ satisfies \emph{(RN4-1):}

Given the solid part of the commutative diagram
$$\xymatrix{A_0 \ar[r]^{a_0}\ar[d]^{f_0} & A_1 \ar[r]^{a_1}\ar[d]^{f_1} & A_2 \ar[r]^{a_2}\ar@{-->}[d]^{f_2} & \cdots \ar[r]^{a_{n-1}}& A_{n} \ar[r]^{a_{n}\hspace{2mm}}\ar@{-->}[d]^{f_{n}} & A_{n+1} \ar[r]^{a_{n+1}}\ar@{-->}[d]^{f_{n+1}} & \Sigma A_0 \ar[d]^{\Sigma f_0}\\
B_0 \ar[r]^{b_0} & B_1 \ar[r]^{b_1} & B_2 \ar[r]^{b_2} & \cdots \ar[r]^{b_{n-1}}& B_{n} \ar[r]^{b_{n}\hspace{2mm}} & B_{n+1} \ar[r]^{b_{n+1}} & \Sigma B_0
}$$
with rows in $\Theta$. Then there exist the dotted morphisms such that the above diagram commutes and the {\bf mapping cone}
$$A_1 \oplus B_0 \xrightarrow{\left[
                    \begin{smallmatrix}
                      -a_1 & 0 \\
                      f_1 & b_0 \\
                    \end{smallmatrix}
                  \right]} A_2\oplus B_1\xrightarrow{\left[
                             \begin{smallmatrix}
                                -a_2 & 0 \\
                                 f_2 & b_1 \\
                             \end{smallmatrix}
                           \right]} A_3\oplus B_2\xrightarrow{\left[
                             \begin{smallmatrix}
                                -a_3 & 0 \\
                                 f_3 & b_2 \\
                             \end{smallmatrix}
                           \right]} \cdots \hspace{50mm}$$
$$\cdots\xrightarrow{\left[
 \begin{smallmatrix}
     -a_n & 0 \\
   f_n & b_{n-1} \\
  \end{smallmatrix}
  \right]} A_{n+1} \oplus B_n \xrightarrow{\left[
                            \begin{smallmatrix}
                                -a_{n+1} & 0 \\
                                 f_{n+1} & b_{n} \\
                             \end{smallmatrix}
                             \right]} \Sigma A_0 \oplus B_{n+1} \xrightarrow{\left[
                            \begin{smallmatrix}
                                -\Sigma a_0 & 0 \\
                                 \Sigma f_0 & b_{n+1} \\
                             \end{smallmatrix}
                             \right]}\Sigma A_1 \oplus \Sigma B_0$$
belongs to $\Theta$.
\end{theorem}
\vspace{-0.6cm}

\proof  First, we  prove (1) implies (2):
Assume that we have the following commutative diagram
$$\xymatrix{A_0 \ar[r]^{a_0}\ar[d]^{f_0} & A_1 \ar[r]^{a_1}\ar[d]^{f_1} & A_2 \ar[r]^{a_2} & \cdots \ar[r]^{a_{n-1}}& A_{n} \ar[r]^{a_{n}\hspace{2mm}} & A_{n+1} \ar[r]^{a_{n+1}} & \Sigma A_0 \ar[d]^{\Sigma f_0}\\
B_0 \ar[r]^{b_0} & B_1 \ar[r]^{b_1} & B_2 \ar[r]^{b_2} & \cdots \ar[r]^{b_{n-1}}& B_{n} \ar[r]^{b_{n}\hspace{2mm}} & B_{n+1} \ar[r]^{b_{n+1}} & \Sigma B_0
}$$
with rows in $\Theta$.
Since $\Theta$ closed under direct sums, then the direct sum of three right $(n+2)$-angles
$$A_0\xrightarrow{~1~} A_0\rightarrow 0\rightarrow \cdots\rightarrow 0\rightarrow \Sigma A_0$$
$$0\rightarrow A_1\xrightarrow{~1~} A_1\rightarrow 0\rightarrow \cdots\rightarrow 0\rightarrow 0$$
$$0\rightarrow B_0\xrightarrow{~1~} B_0\rightarrow 0\rightarrow \cdots\rightarrow 0\rightarrow 0$$
can be combined to  get the following $(n+2)$-$\Sigma$-sequence
$$\xymatrix@C=1.2cm{A_0 \ar[r]^{\left[
                    \begin{smallmatrix}
                      0 \\
                      1 \\
                      0 \\
                    \end{smallmatrix}
                  \right]\hspace{9mm}} & A_1\oplus A_0\oplus B_0 \ar[r]^{\hspace{4mm}\left[
                             \begin{smallmatrix}
                               1 & 0 &0\\
                               0 & 0 &1 \\
                             \end{smallmatrix}
                           \right]} & A_1\oplus B_0 \ar[r] & 0 \ar[r] &\cdots \ar[r] & 0 \ar[r] & \Sigma A_0
}$$
belongs to $\Theta$. By the commutative diagram below
\begin{equation}\label{three sums 1 diagram}
\begin{split}
\xymatrix@C=1.2cm@R=1.7cm{A_0 \ar[r]^{\left[
                    \begin{smallmatrix}
                      0 \\
                      1 \\
                      0 \\
                    \end{smallmatrix}
                  \right]\hspace{9mm}} \ar@{=}[d] & A_1\oplus A_0\oplus B_0 \ar[r]^{\hspace{4mm}\left[
                             \begin{smallmatrix}
                               1 & 0 &0\\
                               0 & 0 &1 \\
                             \end{smallmatrix}
                           \right]} \ar[d]^{\left[
                             \begin{smallmatrix}
                               1 & 0 &0\\
                               0 & -1 &0 \\
                               0 & f_0 &1 \\
                             \end{smallmatrix}
                           \right]} & A_1\oplus B_0 \ar[r] \ar@{=}[d] & 0 \ar[r] \ar@{=}[d] &\cdots \ar[r] & 0 \ar[r] \ar@{=}[d] & \Sigma A_0 \ar@{=}[d]\\
A_0 \ar[r]^{\left[
                    \begin{smallmatrix}
                      0 \\
                      -1 \\
                      f_0 \\
                    \end{smallmatrix}
                  \right]\hspace{9mm}} & A_1\oplus A_0\oplus B_0 \ar[r]^{\hspace{4mm}\left[
                             \begin{smallmatrix}
                               1 & 0 &0\\
                               0 & f_0 &1 \\
                             \end{smallmatrix}
                           \right]} & A_1\oplus B_0 \ar[r] & 0 \ar[r] & \cdots \ar[r] & 0 \ar[r] & \Sigma A_0
}
\end{split}
\end{equation}
we obtained the second row of (\ref{three sums 1 diagram}) belongs to $\Theta$
because $\Theta$ is under isomorphisms.

Similarly, we consider the following commutative diagram
\begin{equation}\label{two sums 2 diagram}
\begin{split}
\xymatrix@C=1.2cm@R=1.7cm{A_0 \ar[r]^{\left[
                    \begin{smallmatrix}
                      a_0 \\
                      0 \\
                    \end{smallmatrix}
                  \right]\hspace{5mm}} \ar@{=}[d] & A_1\oplus B_1 \ar[r]^{\left[
                             \begin{smallmatrix}
                               a_1 &0\\
                               0 &1 \\
                             \end{smallmatrix}
                           \right]} \ar[d]^{\left[
                             \begin{smallmatrix}
                               -1 & 0\\
                               f_1 &1 \\
                             \end{smallmatrix}
                           \right]} & A_2\oplus B_1 \ar[r]^{\hspace{4mm}\left[
                             \begin{smallmatrix}
                               a_2 &~ 0\\
                             \end{smallmatrix}
                           \right]} \ar@{=}[d] & A_3 \ar[r]^{a_3} \ar@{=}[d] &\cdots \ar[r]^{a_n\hspace{2mm}} & A_{n+1} \ar[r]^{a_{n+1}} \ar@{=}[d] & \Sigma A_0 \ar@{=}[d]\\
A_0 \ar[r]^{\left[
                    \begin{smallmatrix}
                      -a_0 \\
                      b_0f_0 \\
                    \end{smallmatrix}
                  \right]\hspace{5mm}} & A_1\oplus B_1 \ar[r]^{\left[
                             \begin{smallmatrix}
                               -a_1 &0\\
                               f_1 &1 \\
                             \end{smallmatrix}
                           \right]} & A_2\oplus B_1 \ar[r]^{\hspace{4mm}\left[
                             \begin{smallmatrix}
                               a_2 &~ 0\\
                             \end{smallmatrix}
                           \right]} & A_3 \ar[r]^{a_3} & \cdots \ar[r]^{a_{n}\hspace{2mm}} & A_{n+1} \ar[r]^{a_{n+1}} & \Sigma A_0
}
\end{split}
\end{equation}
where the first row of (\ref{two sums 2 diagram}) is a direct sum of the following two right $(n+2)$-angles
$$A_0 \xrightarrow{~a_0~} A_1 \xrightarrow{~a_1~} A_2\xrightarrow{~a_2~} \cdots\xrightarrow{~a_n~} A_{n+1}\xrightarrow{~a_{n+1}~} \Sigma A_0,$$
$$0\rightarrow B_1\xrightarrow{~1~} B_1\rightarrow 0\rightarrow \cdots\rightarrow 0\rightarrow 0.$$
Then we get the second row of (\ref{two sums 2 diagram}) belongs to $\Theta$ since $\Theta$ is under isomorphisms.

Consider the commutative diagram below
$$
\xymatrix@C=1.5cm@R=1.7cm{A_1 \oplus A_0\oplus B_0 \ar[r]^{\hspace{4mm}\left[
                    \begin{smallmatrix}
                      1 & 0 & 0 \\
                      0 & 0 & b_0\\
                    \end{smallmatrix}
                  \right]} \ar[d]^{\left[
                    \begin{smallmatrix}
                      1 & -a_0 & 0 \\
                      0 & 1 & 0\\
                      0 & 0 & 1 \\
                    \end{smallmatrix}
                  \right]} & A_1\oplus B_1 \ar[r]^{\hspace{4mm}\left[
                             \begin{smallmatrix}
                              0 &~ b_1\\
                             \end{smallmatrix}
                           \right]} \ar@{=}[d] & B_2 \ar[r]^{b_2} \ar@{=}[d] & \cdots\\
A_1 \oplus A_0\oplus B_0 \ar[r]^{\hspace{4mm}\left[
                    \begin{smallmatrix}
                      1 & a_0 & 0 \\
                      0 & 0 & b_0\\
                    \end{smallmatrix}
                  \right]} & A_1\oplus B_1 \ar[r]^{\hspace{4mm}\left[
                             \begin{smallmatrix}
                              0 &~ b_1\\
                             \end{smallmatrix}
                           \right]} & B_2 \ar[r]^{b_2} & \cdots
}
\hspace{35mm}$$
\begin{equation}\label{three sums 3 diagram}
\begin{split}
\xymatrix@C=2cm@R=1.7cm{\cdots \ar[r]^{b_{n-1}} & B_n \ar[r]^{\left[
                             \begin{smallmatrix}
                              0\\
                               b_n\\
                             \end{smallmatrix}
                           \right]\hspace{6mm}} \ar@{=}[d] & \Sigma A_0\oplus B_{n+1} \ar[r]^{\left[
                             \begin{smallmatrix}
                               0 & 0\\
                               1 & 0\\
                               0 & b_{n+1}\\
                             \end{smallmatrix}
                           \right]\hspace{6mm}} \ar@{=}[d] & \Sigma A_1\oplus \Sigma A_0 \oplus \Sigma B_0 \ar[d]^{\left[
                    \begin{smallmatrix}
                      1 & -\Sigma a_0 & 0 \\
                      0 & 1 & 0\\
                      0 & 0 & 1 \\
                    \end{smallmatrix}
                  \right]}\\
\cdots \ar[r]^{b_{n-1}} & B_n \ar[r]^{\left[
                             \begin{smallmatrix}
                              0\\
                               b_n\\
                             \end{smallmatrix}
                           \right]\hspace{7mm}} & \Sigma A_0\oplus B_{n+1} \ar[r]^{\left[
                             \begin{smallmatrix}
                               -\Sigma a_0 & 0\\
                               1 & 0\\
                               0 & b_{n+1}\\
                             \end{smallmatrix}
                           \right]\hspace{7mm}} & \Sigma A_1\oplus \Sigma A_0 \oplus \Sigma B_0
}
\end{split}
\end{equation}
where the first row of (\ref{three sums 3 diagram}) is a direct sum of the following three right $(n+2)$-angles
$$B_0 \xrightarrow{~b_0~} B_1 \xrightarrow{~b_1~} B_2\xrightarrow{~b_2~} \cdots\xrightarrow{~b_n~} B_{n+1}\xrightarrow{~b_{n+1}~} \Sigma B_0,$$
$$A_1\xrightarrow{1} A_1\rightarrow 0\rightarrow \cdots\rightarrow 0\rightarrow \Sigma A_1,$$
$$A_0\rightarrow 0\rightarrow 0\rightarrow \cdots\rightarrow \Sigma A_0\xrightarrow{1} \Sigma A_0.$$
Then we get the second row of (\ref{three sums 3 diagram}) belongs to $\Theta$ since $\Theta$ is under isomorphisms.
Therefore, we apply (RN4$^*$) to the solid part of the following commutative diagram
$$\xymatrix@R=1.5cm@C=1.3cm{A_0 \ar[r]^{\left[
                    \begin{smallmatrix}
                      0 \\
                      -1 \\
                      f_0 \\
                    \end{smallmatrix}
                  \right]\hspace{8mm}} \ar@{=}[d] & A_1\oplus A_0\oplus B_0 \ar[r]^{\hspace{4mm}\left[
                             \begin{smallmatrix}
                               1 & 0 &0\\
                               0 & f_0 &1 \\
                             \end{smallmatrix}
                           \right]} \ar[d]^{\left[
                             \begin{smallmatrix}
                               1 & a_0 &0\\
                               0 & 0 &b_0 \\
                             \end{smallmatrix}
                           \right]} & A_1\oplus B_0 \ar[r] \ar@{-->}[d]^{u} & 0 \ar[r] \ar@{-->}[d] & \cdots\\
A_0 \ar[r]^{\left[
                    \begin{smallmatrix}
                      -a_0 \\
                      b_0f_0 \\
                    \end{smallmatrix}
                  \right]\hspace{4mm}} \ar[d]^{\left[
                    \begin{smallmatrix}
                      0 \\
                      -1 \\
                      f_0 \\
                    \end{smallmatrix}
                  \right]}& A_1\oplus B_1 \ar[r]^{\left[
                             \begin{smallmatrix}
                               -a_1 &0\\
                               f_1 &1 \\
                             \end{smallmatrix}
                           \right]} \ar@{=}[d] & A_2\oplus B_1 \ar[r]^{\left[
                             \begin{smallmatrix}
                               a_2 & 0\\
                             \end{smallmatrix}
                           \right]} \ar@{-->}[d]^{v_2} & A_3 \ar[r]^{a_3} \ar@{-->}[d]^{v_3} & \cdots \\
A_1 \oplus A_0\oplus B_0 \ar[r]^{\hspace{4mm}\left[
                    \begin{smallmatrix}
                      1 & a_0 & 0 \\
                      0 & 0 & b_0\\
                    \end{smallmatrix}
                  \right]} & A_1\oplus B_1 \ar[r]^{\hspace{4mm}\left[
                             \begin{smallmatrix}
                              0 & b_1\\
                             \end{smallmatrix}
                           \right]} & B_2 \ar[r]^{b_2} & B_3 \ar[r]^{b_3} &\cdots
}\hspace{15mm}$$
\begin{equation}\label{apply rn4^* diagram}
\begin{split}
\hspace{15mm}\xymatrix@R=1.5cm@C=2cm{\cdots \ar[r] & 0 \ar[r] \ar@{-->}[d] & 0 \ar[r] \ar@{-->}[d] & \Sigma A_0 \ar@{=}[d]\\
\cdots \ar[r]^{a_{n-1}} & A_n \ar[r]^{a_n} \ar@{-->}[d]^{v_n} & A_{n+1} \ar[r]^{a_{n+1}} \ar@{-->}[d]^{v_{n+1}} & \Sigma A_0 \ar[d]^{\left[
                    \begin{smallmatrix}
                      0 \\
                      -1 \\
                      \Sigma f_0 \\
                    \end{smallmatrix}
                  \right]}\\
\cdots \ar[r]^{b_{n-1}} & B_n \ar[r]^{\left[
                             \begin{smallmatrix}
                              0\\
                               b_n\\
                             \end{smallmatrix}
                           \right]\hspace{8mm}} & \Sigma A_0\oplus B_{n+1} \ar[r]^{\left[
                             \begin{smallmatrix}
                               -\Sigma a_0 & 0\\
                               1 & 0\\
                               0 & b_{n+1}\\
                             \end{smallmatrix}
                           \right]\hspace{6mm}} & \Sigma A_1\oplus \Sigma A_0 \oplus \Sigma B_0
}
\end{split}
\end{equation}
with rows in $\Theta$. Then there exist morphisms
$$\begin{aligned}
&u: A_1\oplus B_0\rightarrow A_2\oplus B_1 \\&v_2: A_2\oplus B_1\rightarrow B_2\\
&v_i: A_i\rightarrow B_i~~(i=3, \cdots , n) \\&v_{n+1}:A_{n+1}\rightarrow \Sigma A_0\oplus B_{n+1}
\end{aligned}$$
such that (\ref{apply rn4^* diagram}) commutes and the following $(n+2)$-$\Sigma$-sequence
$$A_1 \oplus B_0 \xrightarrow{~u~} A_2\oplus B_1\xrightarrow{\left[
                             \begin{smallmatrix}
                                -a_2 &~ 0 \\
                                 v_{2, 1} &~ v_{2, 2} \\
                             \end{smallmatrix}
                           \right]} A_3\oplus B_2\xrightarrow{\left[
                             \begin{smallmatrix}
                                -a_3 & 0 \\
                                 v_3 & b_2 \\
                             \end{smallmatrix}
                           \right]} \cdots \hspace{50mm}$$
\begin{equation}\label{apply rn4^* diagram make the sequence}
\hspace{10mm}\cdots\xrightarrow{\left[
 \begin{smallmatrix}
     -a_n & 0 \\
   v_n & b_{n-1} \\
  \end{smallmatrix}
  \right]} A_{n+1} \oplus B_n \xrightarrow{\left[
                            \begin{smallmatrix}
                                v_{n+1, 1} &~ 0 \\
                                 v_{n+1, 2} &~ b_{n} \\
                             \end{smallmatrix}
                             \right]} \Sigma A_0 \oplus B_{n+1} \xrightarrow{\left[
                            \begin{smallmatrix}
                                -\Sigma a_0 & 0 \\
                                 \Sigma f_0 & b_{n+1} \\
                             \end{smallmatrix}
                             \right]}\Sigma A_1 \oplus \Sigma B_0
\end{equation}
belongs to $\Theta$.

By the commutative diagram (\ref{apply rn4^* diagram}), we have the following equalities.
$$
\begin{aligned}
\left[
                            \begin{matrix}
                                v_{2, 1} & v_{2, 2} \\
                             \end{matrix}
                             \right]& =\left[
                            \begin{matrix}
                                f_2 & b_1, \\
                             \end{matrix}
                             \right]\\
v_3 & = f_3,\\
&\cdots\\
v_n& =f_n,\\
\left[
                            \begin{matrix}
                                v_{n+1, 1} \\
                                v_{n+1, 2} \\
                             \end{matrix}
                             \right]& =\left[
                            \begin{matrix}
                                -a_{n+1} \\
                                f_{n+1} \\
                             \end{matrix}
                             \right].
\end{aligned}
$$
For $f_i: A_i\rightarrow B_i~(i=2, \cdots, n+1)$,
we have the following commutative diagram, that is to say, $f=(f_0, f_1, f_2, \cdots, f_{n+1})$ is a morphism of right $(n+2)$-angles
$$\xymatrix{A_0 \ar[r]^{a_0}\ar[d]^{f_0} & A_1 \ar[r]^{a_1}\ar[d]^{f_1} & A_2 \ar[r]^{a_2}\ar[d]^{f_2} & \cdots \ar[r]^{a_{n-1}}& A_{n} \ar[r]^{a_{n}\hspace{2mm}}\ar[d]^{f_{n}} & A_{n+1} \ar[r]^{a_{n+1}}\ar[d]^{f_{n+1}} & \Sigma A_0 \ar[d]^{\Sigma f_0}\\
B_0 \ar[r]^{b_0} & B_1 \ar[r]^{b_1} & B_2 \ar[r]^{b_2} & \cdots \ar[r]^{b_{n-1}}& B_{n} \ar[r]^{b_{n}\hspace{2mm}} & B_{n+1} \ar[r]^{b_{n+1}} & \Sigma B_0.
}$$
By
$$u \left[
                             \begin{matrix}
                               1 & 0 &0\\
                               0 & f_0 &1 \\
                             \end{matrix}
                           \right]= \left[
                             \begin{matrix}
                               -a_1 &0\\
                               f_1 &1 \\
                             \end{matrix}
                           \right]\left[
                             \begin{matrix}
                               1 & a_0 &0\\
                               0 & 0 &b_0 \\
                             \end{matrix}
                           \right],$$
we get $u= \left[
                             \begin{matrix}
                               -a_1 &0\\
                               f_1 &b_0 \\
                             \end{matrix}
                           \right]$.
Substituting $u, v_{i,j}$ into (\ref{apply rn4^* diagram make the sequence}), we get
the $(n+2)$-$\Sigma$-sequence
$$A_1 \oplus B_0 \xrightarrow{\left[
                    \begin{smallmatrix}
                      -a_1 & 0 \\
                      f_1 & b_0 \\
                    \end{smallmatrix}
                  \right]} A_2\oplus B_1\xrightarrow{\left[
                             \begin{smallmatrix}
                                -a_2 & 0 \\
                                 f_2 & b_1 \\
                             \end{smallmatrix}
                           \right]} A_3\oplus B_2\xrightarrow{\left[
                             \begin{smallmatrix}
                                -a_3 & 0 \\
                                 f_3 & b_2 \\
                             \end{smallmatrix}
                           \right]} \cdots \hspace{50mm}$$
$$\cdots\xrightarrow{\left[
 \begin{smallmatrix}
     -a_n & 0 \\
   f_n & b_{n-1} \\
  \end{smallmatrix}
  \right]} A_{n+1} \oplus B_n \xrightarrow{\left[
                            \begin{smallmatrix}
                                -a_{n+1} & 0 \\
                                 f_{n+1} & b_{n} \\
                             \end{smallmatrix}
                             \right]} \Sigma A_0 \oplus B_{n+1} \xrightarrow{\left[
                            \begin{smallmatrix}
                                -\Sigma a_0 & 0 \\
                                 \Sigma f_0 & b_{n+1} \\
                             \end{smallmatrix}
                             \right]}\Sigma A_1 \oplus \Sigma B_0$$
belongs to $\Theta$.

Next, we prove (2) implies (1):
This proof process is similar to the necessity of proving Theorem \ref{th N4 equivalent}, and we have omitted it here. \qed

\begin{remark}
The above result
 was proved in \cite[Thoerem 4.4]{bt} for an $(n+2)$-angulated category.
Now it has been extended to a right  $(n+2)$-angulated category.
However, our proof method differs from \cite[Thoerem 4.4]{bt} as $\Sigma$ is not automorphism for
a right $(n+2)$-angulated category.
\end{remark}

\begin{theorem}\label{main6}
If $\Theta$ is a collection of $(n+2)$-$\Sigma$-sequences satisfying the axioms \emph{(RN1)}, {\rm (RN2)} and {\rm (RN3)}, then the following statements are equivalent:

$(1)$ $\Theta$ satisfies \emph{(RN4-1);}

$(2)$ $\Theta$ satisfies \emph{(RN4-2):}

Given the solid part of the commutative diagram
$$
\xymatrix{A_0 \ar[r]^{a_0}\ar@{=}[d] & A_1 \ar[r]^{a_1}\ar[d]^{f_1} & A_2 \ar[r]^{a_2} \ar@{-->}[d]^{f_2} & \cdots \ar[r]^{a_{n-1}}& A_{n} \ar[r]^{a_{n}\hspace{2mm}} \ar@{-->}[d]^{f_n} & A_{n+1} \ar[r]^{a_{n+1}} \ar@{-->}[d]^{f_{n+1}} & \Sigma A_0 \ar@{=}[d]\\
A_0 \ar[r]^{b_0} & B_1 \ar[r]^{b_1} & B_2 \ar[r]^{b_2} & \cdots \ar[r]^{b_{n-1}}& B_{n} \ar[r]^{b_{n}\hspace{2mm}} & B_{n+1} \ar[r]^{b_{n+1}} & \Sigma A_0
}
$$
with rows in $\Theta$. Then there exist the dotted morphisms $f_i: A_i\rightarrow B_i~(i=2, \cdots, n+1)$ such that the above diagram commutes and the following $(n+2)$-$\Sigma$-sequence
$$A_1 \xrightarrow{\left[
                    \begin{smallmatrix}
                      -a_1 \\
                      f_1 \\
                    \end{smallmatrix}
                  \right]} A_2\oplus B_1\xrightarrow{\left[
                             \begin{smallmatrix}
                                -a_2 & 0 \\
                                 f_2 & b_1 \\
                             \end{smallmatrix}
                           \right]} A_3\oplus B_2\xrightarrow{\left[
                             \begin{smallmatrix}
                                -a_3 & 0 \\
                                 f_3 & b_2 \\
                             \end{smallmatrix}
                           \right]} \cdots \hspace{50mm}$$
$$\cdots\xrightarrow{\left[
 \begin{smallmatrix}
     -a_n & 0 \\
      f_n & b_{n-1} \\
  \end{smallmatrix}
  \right]} A_{n+1} \oplus B_n \xrightarrow{\left[
                            \begin{smallmatrix}
                                 f_{n+1} &~ b_{n} \\
                             \end{smallmatrix}
                             \right]} B_{n+1} \xrightarrow{\Sigma a_0\circ b_{n+1}}\Sigma A_1 $$
belongs to $\Theta$.
\end{theorem}

\proof First, we prove (1) implies (2): Assume that we have the solid part of the following commutative
\begin{equation}\label{rn4-2 orignal diagram}
\begin{split}
\xymatrix{A_0 \ar[r]^{a_0}\ar@{=}[d] & A_1 \ar[r]^{a_1}\ar[d]^{f_1} & A_2 \ar[r]^{a_2} \ar@{-->}[d]^{f_2} & \cdots \ar[r]^{a_{n-1}}& A_{n} \ar[r]^{a_{n}\hspace{2mm}} \ar@{-->}[d]^{f_n} & A_{n+1} \ar[r]^{a_{n+1}} \ar@{-->}[d]^{f_{n+1}} & \Sigma A_0 \ar@{=}[d]\\
A_0 \ar[r]^{b_0} & B_1 \ar[r]^{b_1} & B_2 \ar[r]^{b_2} & \cdots \ar[r]^{b_{n-1}}& B_{n} \ar[r]^{b_{n}\hspace{2mm}} & B_{n+1} \ar[r]^{b_{n+1}} & \Sigma A_0
}
\end{split}
\end{equation}
with rows in $\Theta$. Apply (RN4-1) for diagram (\ref{rn4-2 orignal diagram}), there exist $f_i: A_i\rightarrow B_i~(i=2, \cdots, n+1)$ such that (\ref{rn4-2 orignal diagram}) commutes and the mapping cone
$$A_1 \oplus B_0 \xrightarrow{\left[
                    \begin{smallmatrix}
                      -a_1 & 0 \\
                      f_1 & b_0 \\
                    \end{smallmatrix}
                  \right]} A_2\oplus B_1\xrightarrow{\left[
                             \begin{smallmatrix}
                                -a_2 & 0 \\
                                 f_2 & b_1 \\
                             \end{smallmatrix}
                           \right]} A_3\oplus B_2\xrightarrow{\left[
                             \begin{smallmatrix}
                                -a_3 & 0 \\
                                 f_3 & b_2 \\
                             \end{smallmatrix}
                           \right]} \cdots \hspace{50mm}$$
$$\cdots\xrightarrow{\left[
 \begin{smallmatrix}
     -a_n & 0 \\
   f_n & b_{n-1} \\
  \end{smallmatrix}
  \right]} A_{n+1} \oplus B_n \xrightarrow{\left[
                            \begin{smallmatrix}
                                -a_{n+1} & 0 \\
                                 f_{n+1} & b_{n} \\
                             \end{smallmatrix}
                             \right]} \Sigma A_0 \oplus B_{n+1} \xrightarrow{\left[
                            \begin{smallmatrix}
                                -\Sigma a_0 & 0 \\
                                 1 & b_{n+1} \\
                             \end{smallmatrix}
                             \right]}\Sigma A_1 \oplus \Sigma B_0$$
belongs to $\Theta$.

Since $\Theta$ is closed under direct summands, then the following commutative diagram  \newpage
$$\xymatrix@C=2cm@R=1.5cm{A_1 \ar[r]^{\left[
                    \begin{smallmatrix}
                      -a_1 \\
                      f_1 \\
                    \end{smallmatrix}
                  \right]\hspace{4mm}} \ar[d]^{\left[
                    \begin{smallmatrix}
                      1 \\
                      0 \\
                    \end{smallmatrix}
                  \right]} & A_2\oplus B_1 \ar[r]^{\hspace{4mm}\left[
                             \begin{smallmatrix}
                                -a_2 & 0 \\
                                 f_2 & b_1 \\
                             \end{smallmatrix}
                           \right]\hspace{4mm}} \ar@{=}[d] & A_3\oplus B_2 \ar[r]^{\hspace{4mm}\left[
                             \begin{smallmatrix}
                                -a_3 & 0 \\
                                 f_3 & b_2 \\
                             \end{smallmatrix}
                           \right]} \ar@{=}[d] & \cdots\\
A_1\oplus A_0 \ar[r]^{\left[
                    \begin{smallmatrix}
                      -a_1 & 0 \\
                      f_1 & b_0 \\
                    \end{smallmatrix}
                  \right]} \ar[d]^{\left[
                    \begin{smallmatrix}
                      1 & a_0 \\
                    \end{smallmatrix}
                  \right]}& A_2\oplus B_1 \ar[r]^{\hspace{4mm}\left[
                             \begin{smallmatrix}
                                -a_2 & 0 \\
                                 f_2 & b_1 \\
                             \end{smallmatrix}
                           \right]\hspace{4mm}} \ar@{=}[d] & A_3\oplus B_2 \ar[r]^{\hspace{4mm}\left[
                             \begin{smallmatrix}
                                -a_3 & 0 \\
                                 f_3 & b_2 \\
                             \end{smallmatrix}
                           \right]} \ar@{=}[d] &\cdots\\
A_1 \ar[r]^{\left[
                             \begin{smallmatrix}
                                -a_1 \\
                                 f_1 \\
                             \end{smallmatrix}
                           \right]\hspace{4mm}} & A_2\oplus B_1 \ar[r]^{\left[
                             \begin{smallmatrix}
                                -a_2 & 0 \\
                                 f_2 & b_1 \\
                             \end{smallmatrix}
                           \right]} & A_3\oplus B_2 \ar[r]^{\hspace{4mm}\left[
                             \begin{smallmatrix}
                                -a_3 & 0 \\
                                 f_3 & b_2 \\
                             \end{smallmatrix}
                           \right]}&\cdots
}\hspace{25mm}$$
\begin{equation}\label{mapping cone direct sum diagram 6}
\begin{split}
\xymatrix@C=2cm@R=1.5cm{\cdots \ar[r]^{\left[
                             \begin{smallmatrix}
                                -a_n & 0 \\
                                 f_n & b_{n-1} \\
                             \end{smallmatrix}
                           \right]\hspace{7mm}} & A_{n+1}\oplus B_n \ar[r]^{\hspace{7mm}\left[
                            \begin{smallmatrix}
                                 f_{n+1} &~ b_{n} \\
                             \end{smallmatrix}
                             \right]\hspace{4mm}} \ar@{=}[d] & B_{n+1} \ar[r]^{\Sigma a_0 \circ b_{n+1}} \ar[d]^{\left[
                            \begin{smallmatrix}
                                 -b_{n+1} \\
                                 1 \\
                             \end{smallmatrix}
                             \right]} & \Sigma A_1 \ar[d]^{\left[
                            \begin{smallmatrix}
                                 1 \\
                                 0 \\
                             \end{smallmatrix}
                             \right]}\\
 \cdots \ar[r]^{\left[
                             \begin{smallmatrix}
                                -a_n & 0 \\
                                 f_n & b_{n-1} \\
                             \end{smallmatrix}
                           \right]\hspace{7mm}} & A_{n+1}\oplus B_n \ar[r]^{\left[
                             \begin{smallmatrix}
                                -a_{n+1} & 0 \\
                                 f_{n+1} & b_{n} \\
                             \end{smallmatrix}
                           \right]\hspace{2mm}} \ar@{=}[d] & \Sigma A_0\oplus B_{n+1} \ar[r]^{\hspace{6mm}\left[
                            \begin{smallmatrix}
                                 -\Sigma a_0 & 0 \\
                                 1 & b_{n+1}\\
                             \end{smallmatrix}
                             \right]\hspace{6mm}} \ar[d]^{\left[
                            \begin{smallmatrix}
                                 0 & 1\\
                             \end{smallmatrix}
                             \right]} & \Sigma A_1\oplus \Sigma A_0 \ar[d]^{\left[
                            \begin{smallmatrix}
                                 1 & \Sigma a_0\\
                             \end{smallmatrix}
                             \right]}\\
\cdots \ar[r]^{\left[
                             \begin{smallmatrix}
                                -a_n & 0 \\
                                 f_n & b_{n-1} \\
                             \end{smallmatrix}
                           \right]\hspace{7mm}} & A_{n+1}\oplus B_n \ar[r]^{\hspace{7mm}\left[
                            \begin{smallmatrix}
                                 f_{n+1} &~ b_{n} \\
                             \end{smallmatrix}
                             \right]\hspace{4mm}} & B_{n+1} \ar[r]^{\Sigma a_0\circ b_{n+1}} & \Sigma A_1
}
\end{split}
\end{equation}
shows that the first row of (\ref{mapping cone direct sum diagram 6}) belongs to $\Theta$.

Next, we prove (2) implies (1):  Assume that we have the following commutative diagram
\begin{equation}\label{th6.2 2}
\begin{split}
\xymatrix{A_0 \ar[r]^{a_0}\ar[d]^{f_0} & A_1 \ar[r]^{a_1}\ar[d]^{f_1} & A_2 \ar[r]^{a_2} & \cdots \ar[r]^{a_{n-1}}& A_{n} \ar[r]^{a_{n}\hspace{2mm}} & A_{n+1} \ar[r]^{a_{n+1}} & \Sigma A_0 \ar[d]^{\Sigma f_0}\\
B_0 \ar[r]^{b_0} & B_1 \ar[r]^{b_1} & B_2 \ar[r]^{b_2} & \cdots \ar[r]^{b_{n-1}}& B_{n} \ar[r]^{b_{n}\hspace{2mm}} & B_{n+1} \ar[r]^{b_{n+1}} & \Sigma B_0
}
\end{split}
\end{equation}
with rows in $\Theta$.

The $(n+2)$-$\Sigma$-sequence
$$A_0\oplus B_0 \xrightarrow{\left[
                            \begin{smallmatrix}
                                 0 &~ b_0\\
                             \end{smallmatrix}
                             \right]} B_1 \xrightarrow{~b_1~} B_2 \xrightarrow{~b_2~} \cdots \xrightarrow{~b_{n-1}~} B_n \xrightarrow{\left[
                            \begin{smallmatrix}
                                 0 \\
                                 b_n \\
                             \end{smallmatrix}
                             \right]}\Sigma A_0\oplus B_{n+1} \xrightarrow{\left[
                            \begin{smallmatrix}
                                 -1 & 0\\
                                 0 & b_{n+1}\\
                             \end{smallmatrix}
                             \right]}\Sigma A_0\oplus \Sigma B_0$$
belongs to $\Theta$ since it is the direct sum of a right $(n+2)$-angle
$$A_0\rightarrow 0\rightarrow \cdots \rightarrow \Sigma A_0\xrightarrow{-1}\Sigma A_0$$
and the second row of diagram (\ref{th6.2 2})
$$\xymatrix{B_0 \ar[r]^{b_0} & B_1 \ar[r]^{b_1} & B_2 \ar[r]^{b_2} & \cdots \ar[r]^{b_{n-1}}& B_{n} \ar[r]^{b_{n}\hspace{2mm}} & B_{n+1} \ar[r]^{b_{n+1}} & \Sigma B_0}.$$

By the commutative diagram below  \newpage
$$\xymatrix@C=1.5cm@R=1.3cm{A_0\oplus B_0 \ar[r]^{\hspace{4mm}\left[
                            \begin{smallmatrix}
                                 b_0f_0 & b_0\\
                             \end{smallmatrix}
                             \right]}\ar[d]^{\left[
                            \begin{smallmatrix}
                                 1 & 0\\
                                 f_0 & 1\\
                             \end{smallmatrix}
                             \right]} & B_1 \ar[r]^{b_1} \ar@{=}[d] & B_2 \ar[r]^{b_2} \ar@{=}[d] & \cdots \\
A_0\oplus B_0 \ar[r]^{\hspace{4mm}\left[
                            \begin{smallmatrix}
                                 0 & b_0\\
                             \end{smallmatrix}
                             \right]} & B_1 \ar[r]^{b_1} & B_2 \ar[r]^{b_2} & \cdots
}\hspace{50mm}
$$
\begin{equation}\label{th6.2 isomorphism diagram}
\begin{split}
\hspace{20mm}\xymatrix@C=1.5cm@R=1.3cm{\cdots \ar[r]^{b_{n-1}}& B_{n} \ar[r]^{\left[
                            \begin{smallmatrix}
                                 0 \\
                                 b_n \\
                             \end{smallmatrix}
                             \right]\hspace{6mm}} \ar@{=}[d] & \Sigma A_0\oplus B_{n+1} \ar[r]^{\left[
                            \begin{smallmatrix}
                                 -1 & 0\\
                                 \Sigma f_0 & b_{n+1}\\
                             \end{smallmatrix}
                             \right]} \ar@{=}[d] & \Sigma A_0\oplus \Sigma B_0 \ar[d]^{\left[
                            \begin{smallmatrix}
                                 1 & 0\\
                                \Sigma f_0 & 1\\
                             \end{smallmatrix}
                             \right]}\\
\cdots \ar[r]^{b_{n-1}}& B_{n} \ar[r]^{\left[
                            \begin{smallmatrix}
                                 0 \\
                                 b_n \\
                             \end{smallmatrix}
                             \right]\hspace{6mm}} & \Sigma A_0\oplus B_{n+1} \ar[r]^{\left[
                            \begin{smallmatrix}
                                 -1 & 0\\
                                 0 & b_{n+1}\\
                             \end{smallmatrix}
                             \right]} & \Sigma A_0\oplus \Sigma B_0,
}
\end{split}
\end{equation}
we get the first row of (\ref{th6.2 isomorphism diagram}) belongs to $\Theta$.

Consider the following commutative diagram
$$\xymatrix@C=1.5cm@R=1.3cm{A_0\oplus B_0 \ar[r]^{\left[
                            \begin{smallmatrix}
                                 a_0 & 0\\
                                 0 & 1 \\
                             \end{smallmatrix}
                             \right]}\ar@{=}[d] & A_1\oplus B_0 \ar[r]^{\hspace{4mm}\left[
                            \begin{smallmatrix}
                                 a_1 & 0\\
                             \end{smallmatrix}
                             \right]} \ar[d]^{\left[
                            \begin{smallmatrix}
                                 f_1 & b_0\\
                             \end{smallmatrix}
                             \right]} & A_2 \ar[r]^{a_2} & \cdots \\
A_0\oplus B_0 \ar[r]^{\hspace{4mm}\left[
                            \begin{smallmatrix}
                                 b_0f_0 & b_0\\
                             \end{smallmatrix}
                             \right]} & B_1 \ar[r]^{b_1} & B_2 \ar[r]^{b_2} & \cdots
}\hspace{50mm}
$$
\begin{equation}\label{th6.2 apply rn4-2 6.19}
\begin{split}
\hspace{20mm}\xymatrix@C=1.5cm@R=1.3cm{\cdots \ar[r]^{a_{n-1}}& A_{n} \ar[r]^{a_n} & A_{n+1} \ar[r]^{\left[
                            \begin{smallmatrix}
                                 a_{n+1}\\
                                 0\\
                             \end{smallmatrix}
                             \right]\hspace{6mm}} & \Sigma A_0\oplus \Sigma B_0 \ar@{=}[d]\\
\cdots \ar[r]^{b_{n-1}}& B_{n} \ar[r]^{\left[
                            \begin{smallmatrix}
                                 0 \\
                                 b_n \\
                             \end{smallmatrix}
                             \right]\hspace{8mm}} & \Sigma A_0\oplus B_{n+1} \ar[r]^{\left[
                            \begin{smallmatrix}
                                 -1 & 0\\
                                 \Sigma f_0 & b_{n+1}\\
                             \end{smallmatrix}
                             \right]} & \Sigma A_0\oplus \Sigma B_0
}
\end{split}
\end{equation}
where the first row of (\ref{th6.2 apply rn4-2 6.19}) obtained by the direct sum of
$$B_0\xrightarrow{1} B_0\rightarrow 0\rightarrow \cdots \rightarrow 0\rightarrow\Sigma B_0$$
and the first row of diagram (\ref{th6.2 2}), we
apply (RN4-2) for (\ref{th6.2 apply rn4-2 6.19}), there exist
$$f_i: A_i\rightarrow B_i~(i=2, \cdots, n), ~~\alpha: A_{n+1}\rightarrow \Sigma A_0\oplus B_{n+1}$$ such that the diagram (\ref{th6.2 apply rn4-2 6.19}) commutes, where $\alpha=\left[
                            \begin{matrix}
                                 -a_{n+1} \\
                                 f_{n+1} \\
                             \end{matrix}
                             \right]$. i.e. $f=(f_0, f_1, f_2, \cdots, f_{n+1})$ is a morphism of right $(n+2)$-angles in (\ref{th6.2 2}),
and the the $(n+2)$-$\Sigma$-sequence
$$A_1 \oplus B_0 \xrightarrow{\left[
                    \begin{smallmatrix}
                      -a_1 & 0 \\
                      f_1 & b_0 \\
                    \end{smallmatrix}
                  \right]} A_2\oplus B_1\xrightarrow{\left[
                             \begin{smallmatrix}
                                -a_2 & 0 \\
                                 f_2 & b_1 \\
                             \end{smallmatrix}
                           \right]} A_3\oplus B_2\xrightarrow{\left[
                             \begin{smallmatrix}
                                -a_3 & 0 \\
                                 f_3 & b_2 \\
                             \end{smallmatrix}
                           \right]} \cdots \hspace{50mm}$$
$$\cdots\xrightarrow{\left[
 \begin{smallmatrix}
     -a_n & 0 \\
   f_n & b_{n-1} \\
  \end{smallmatrix}
  \right]} A_{n+1} \oplus B_n \xrightarrow{\left[
                            \begin{smallmatrix}
                                -a_{n+1} & 0 \\
                                 f_{n+1} & b_{n} \\
                             \end{smallmatrix}
                             \right]} \Sigma A_0 \oplus B_{n+1} \xrightarrow{\left[
                            \begin{smallmatrix}
                                -\Sigma a_0 & 0 \\
                                 \Sigma f_0 & b_{n+1} \\
                             \end{smallmatrix}
                             \right]}\Sigma A_1 \oplus \Sigma B_0$$
belongs to $\Theta$. This completes the proof.\qed

\begin{remark}
Note that Theorem \ref{main6} was proved in \cite[Thoerem 3.1]{LZ19} for an $(n+2)$-angulated category.
Now it has been extended to a right $(n+2)$-angulated category.
However, our proof method differs from \cite[Thoerem 3.1]{LZ19} as $\Sigma$ is not automorphism for
a right $(n+2)$-angulated category.
\end{remark}

%\hspace{-4mm}\textbf{Data Availability}\hspace{2mm} Data sharing not applicable to this article as no datasets were generated or analysed during
%the current study.
%\vspace{2mm}
%
%\hspace{-4mm}\textbf{Conflict of Interests}\hspace{2mm} The authors declare that they have no conflicts of interest to this work.

\textbf{Jing He}\\
School of Science, Hunan University of Technology and Business, 410205 Changsha, Hunan, P. R. China\\
E-mail: \textsf{jinghe1003@163.com}\\[0.3cm]
\textbf{Jiangsha Li}\\
College of Mathematics, Hunan Institute of Science and Technology, 414006 Yueyang, Hunan, P. R. China.\\
E-mail: \textsf{jiangshali322@163.com}

\end{document}